\documentclass[12pt]{article}

\usepackage{amsmath}
\usepackage{amssymb}
\usepackage{amsthm}
\usepackage{latexsym}
\usepackage{color}
\usepackage{graphicx}
\usepackage{color}

\DeclareSymbolFont{calletters}{OMS}{cmsy}{m}{n}
\DeclareSymbolFontAlphabet{\mathcal}{calletters}

\def\be{\begin{eqnarray}}
\def\ee{\end{eqnarray}}

\def\b*{\begin{eqnarray*}}
\def\e*{\end{eqnarray*}}

\newtheorem{Theorem}{Theorem}[section]

\newtheorem{Proposition}[Theorem]{Proposition}

\newtheorem{Assumption}[Theorem]{Assumption}

\newtheorem{Remark}[Theorem]{Remark}
\newtheorem{Example}[Theorem]{Example}

\makeatletter \@addtoreset{equation}{section}

\newcommand{\un}{1\hspace{-1mm}{\rm I}}   

\usepackage{color}

\newcommand{\cblue}{\color{blue}}

\newcommand{\cred}{\color{red}}

\def \E{\mathbb{E}}

\def \L{\mathbb{L}}
\def \P{\mathbb{P}}

\def \R{\mathbb{R}}
\def \S{\mathbb{S}}
\def \M{\mathbb{M}}
\def \N{\mathbb{N}}

\def\Fc{{\cal F}}

\def\Kc{{\cal K}}
\def\Lc{{\cal L}}

\def\Wc{{\cal W}}
\def\Xc{{\cal X}}

\def\Lb{\overline L}

\def \I{{\bf I}}

\def \Om{\Omega}

\def \eps{\varepsilon}

\def \0{\mathbf{0}}

\def \Xb{\overline{X}}

\def \Kcb{\overline{\Kc}}
\def \psih{\widehat{\psi}}

\newcommand{\rmi}{{\rm (i)$\>\>$}}

\newcommand{\rmii}{{\rm (ii)$\>\>$}}
\newcommand{\rmiii}{{\rm (iii)$\>\,    \,$}}
\newcommand{\rmiv}{{\rm (iv)$\>\>$}}

\def\x{\times}

\def\1{{\bf 1}}
\def \proof{{\noindent \bf Proof. }}

\def \Wch{\widehat{\Wc}}
\def \Wcb{\overline{\Wc}}
\def \Fb{\overline F}
\def \Lh{\widehat L}
\def \ph{\hat p}
\def \Ih{\widehat I}
\def \Xh{\widehat X}
\def \Wh{\widehat W}
\def \Kch{\widehat \Kc}
\def \Kcbh{\widehat \Kcb}
\def \Th{\widehat T}

\title{Branching diffusion representation of semilinear PDEs and Monte Carlo approximation
\footnote{We are grateful to Vincent Bansaye, Julien Claisse, Emmanuel Gobet and Gaoyue Guo for valuable comments and suggestions.
X. Tan and N. Touzi gratefully acknowledge the financial support of the ERC 321111 Rofirm, the ANR Isotace, and the Chairs Financial Risks (Risk Foundation, sponsored by Soci\'et\'e G\'en\'erale) and Finance and Sustainable Development (IEF sponsored by EDF and CA).}}

\author{Pierre Henry-Labord\`ere\thanks{Soci\'et\'e G\'en\'erale, Global Market Quantitative Research,
pierre.henry-labordere@sgcib.com}
	\and Nadia Oudjane \thanks{EDF R\&D \& FiME, Laboratoire de Finance des March\'es de l'Energie, nadia.oudjane@edf.fr}
	\and Xiaolu Tan\thanks{CEREMADE, University of Paris-Dauphine, PSL Research University, tan@ceremade.dauphine.fr}
	\and Nizar Touzi\thanks{Ecole Polytechnique Paris, Centre de Math\'ematiques Appliqu\'ees, nizar.touzi@polytechnique.edu}
	\and Xavier Warin \thanks{EDF R\&D \& FiME, Laboratoire de Finance des March\'es de l'Energie, xavier.warin@edf.fr}
}

\date{\today}

\begin{document}

\maketitle

\abstract{
	We provide a representation result of parabolic semi-linear PD-Es, with polynomial nonlinearity, by branching diffusion processes.
	We extend the classical representation for KPP equations, introduced by
	Skorokhod \cite{Skorokhod}, Watanabe \cite{Watanabe} and McKean \cite{McKean_1975},
	by allowing for polynomial nonlinearity in the pair $(u, Du)$, where $u$ is the solution of the PDE with space gradient $Du$.
	Similar to the previous literature, our result requires a non-explosion condition which restrict to ``small maturity'' or ``small nonlinearity'' of the PDE.
	Our main ingredient is the automatic differentiation technique as in \cite{HTT2}, based on the Malliavin integration by parts, which allows to account for the nonlinearities in the gradient. As a consequence, the particles of our branching diffusion are marked by the nature of the nonlinearity. This new representation has very important numerical implications as it is suitable for Monte Carlo simulation. Indeed, this provides the first numerical method for high dimensional nonlinear PDEs with error estimate induced by the dimension-free Central limit theorem. The complexity is also easily seen to be of the order of the squared dimension. The final section of this paper illustrates the efficiency of the algorithm by some high dimensional numerical experiments.

\vspace{5mm}

 }

\vspace{2mm}

\section{Introduction}

The objective of the present paper is to provide a probabilistic representation for the solution of a nonlinear parabolic second order partial differential equation (PDE) which is suitable for a high dimensional Monte Carlo approximating scheme. Our main results achieve this goal in the context of semilinear PDEs:
 \b*
 -\partial_tu-\Lc u =f(u,Du), &u_T=g,& t<T,~x\in\R^d,
 \e*
with polynomial non-linearity $f_{t,x}(y,z)$ in the solution and its gradient, diffusion generator $\Lc$, and bounded terminal condition $g$.

Previous representation results were obtained in the literature by means of backward stochastic differential equations, as introduced by Pardoux and Peng \cite{PardouxPeng}. The Monte Carlo numerical implications of this representation were introduced by Bally \& Pag\`es \cite{BallyPages}, Bouchard \& Touzi \cite{BouchardTouzi} and Zhang \cite{Zhang}, and generated a large stream of the literature. However, these methods can be viewed as a Monte Carlo version of the finite elements methods, and as such, are subject to the problem of curse of dimensionality. Our primary goal is to avoid this numerical problem so as to be capable to handle genuinely high-dimensional problems. This however will be achieved at the cost of some limitations...

Our main representation result is obtained by using the branching diffusion trick to absorb the nonlinearity, as illustrated by Skorokhod \cite{Skorokhod}, Watanabe \cite{Watanabe} and McKean \cite{McKean_1975} in the context of the KPP equation, see also the extensions in Rasulov, Raimova \& Mascagni \cite{RasulovRaimovaMascagni} and our previous paper \cite{HTT1} where the representation is also shown to allow for path-dependency.

Since the gradient is also involved in the nonlinearity, our representation result is a significant improvement of the classically well-know representation of KPP equations. We observe that the polynomial nonlinearity naturally induces some restrictions needed to ensure the non-explosion of the corresponding solution. As a consequence, our representation holds under technical conditions of small maturity or small nonlinearity of the PDE.

The main idea for our representation is to use the Monte Carlo automatic differentiation technique in addition to the branching diffusion representation. The automatic differentiation in Monte Carlo approximation of diffusions was successfully used in the previous literature by Fourni\'e et al. \cite{FLLLT}, Bouchard, Ekeland \& Touzi \cite{BET}, Henry-Labord\`ere, Tan \& Touzi \cite{HTT2}, and Doumbia, Oudjane \& Warin \cite{DOW}. The resulting branching diffusion in the representation differs from that of the original founding papers \cite{Skorokhod, Watanabe, McKean_1975} by introducing marks for the particles born at each branching. The mark of the particle determines the nature of the differentiation, and thus induces the corresponding automatic differentiation weight.

We next illustrate the main idea behind our representation in the context of the following extension of the one-dimensional Burgers equation:
 \b*
 \Lc u:=\frac12\Delta u,
 &\mbox{and}&
 f(y,z)=\frac12(y^2+yz).
 \e*
Let $W^1$ be a Brownian motion, and $\tau^1$ an independent random variable with density $\rho>0$ on $\R_+$, and denote $\bar F(t):=\int_t^\infty\rho(s)ds$. We also introduce another independent random variable $I^1$ which takes the values $0$ and $1$ with equal probability. Then, denoting by $\E_{t,x}$ the expectation operator conditional on the starting data $W_t=x$ at time $t$, we obtain from the Feynman-Kac formula the representation of the solution $u$ as:
 $$
 u(0,x)
 =
 \E_{0,x} \Big[\bar F(T)\frac{g(W_T)}{\bar F(T)}+\int_0^T \frac{f(u,Du)(t,W_t)}{\rho(t)}\rho(t)dt\Big]
 =
 \E_{0,x} \big[\phi\big(T_{(1)},W^1_{T_{(1)}}\big)\big],
  $$
where $T_{(1)}:=\tau^1\wedge T$, and
  \be
  \phi(t,y)
  &:=& \frac{\1_{\{t\ge T\}}}{\bar F(T)}g(y)
                    \!+\! \frac{\1_{\{t<T\}}}{\rho(t)}(uD^{I_1}u)(t,y).
 \label{phi}
 \ee
We next consider the two alternative cases for the value of $I^1$.
\begin{itemize}
\item On the event set $\{I^1=0\}$, it follows from the Markov property that:
 $$
 (uD^{I^1}u)(t,y)
 =
 u(t,y)^2
 =
 \E_{t,y}\big[\phi(t+\tau^1,W^1_{t+\tau^1})\big]^2.
 $$
The tricky branching diffusion representation now pops up naturally by rewriting the last expression in terms of independent copies $(W^{1,1},\tau^{1,1})$ and $(W^{1,2},\tau^{1,2})$ as:
 \b*
 (uD^{I_1}u)(t,y)
 &=&
 \E_{t,y}\big[\phi\big(t+\tau^{1,1},W^{1,1}_{t+\tau^{1,1}}\big)\big]
 \E_{t,y}\big[\phi\big(t+\tau^{1,2},W^{1,1}_{t+\tau^{1,2}}\big)\big]
 \\
 &=&
 \E_{t,y}\big[\phi\big(t+\tau^{1,1},W^{1,1}_{t+\tau^{1,1}}\big)\phi\big(t+\tau^{1,2},W^{1,2}_{t+\tau^{1,2}}\big)\big].
 \e*
where $\E_{t,y}$ denotes the expectation operator conditional on $W^{1,1}_t=W^{1,2}_t=y$. Substituting this expression in \eqref{phi} and using the tower property, we see that the branching mechanism allows to absorb the nonlinearity.
\item On the event set $\{I^1=1\}$, we arrive similarly to the expression
  \b*
 (uD^{I_1}u)(t,y)
 &=&
 \E_{t,y}\big[\phi\big(t+\tau^{1,1},W^{1,1}_{t+\tau^{1,1}}\big)\big]
 \partial_y\E_{t,y}\big[\phi\big(t+\tau^{1,2},W^{1,2}_{t+\tau^{1,2}}\big)\big].
 \e*
Our main representation is based on the following automatic differentiation:
 \b*
 \partial_y\E_{t,y}\big[\phi\big(t+\tau^{1,2},W^{1,2}_{t+\tau^{1,2}}\big)\big]
 &=&
 \E_{t,y}\Big[\frac{W^{1,2}_{t+\tau^{1,2}}-W^{1,2}_t}
{\tau^{1,2}}\phi\big(t+\tau^{1,2},W^{1,1}_{t+\tau^{1,2}}\big)\Big],
 \e*
which is an immediate consequence of the differentiation with respect to the heat kernel, i.e. the marginal density of the Brownian motion. By the independence of $W^{1,1}$ and $W^{1,2}$, this provides:
\b*
 (uD^{I_1}u)(t,y)
 &=&
 \E_{t,y}\Big[\frac{W^{1,2}_{t+\tau^{1,2}}-W^{1,2}_t}
{\tau^{1,2}}\phi\big(t+\tau^{1,1},W^{1,1}_{t+\tau^{1,1}}\big)\phi\big(t+\tau^{1,2},W^{1,2}_{t+\tau^{1,2}}\big)\Big],
 \e*
so that the branching mechanism allows again to absorb the nonlinearity by substituting in \eqref{phi} and using the tower property.
\end{itemize}
The two previous cases are covered by denoting $T_{(1,i)}:=T\wedge(\tau^1+\tau^{1,i})$ for $i=0,1$, and introducing the random variable:
 $$
 \Wc^{1}
 \!:=\!
 \1_{\{I^1=0\}} + \1_{\{I^1=1\}}\frac{\Delta W^{1,2}_{T_{(1,2)}}}{\Delta T_{(1,2)}},
 ~\mbox{with}~
 \Delta W^{1,2}_{T_{(1,2)}}\!:=\!W^{1,2}_{T_{(1,2)}}-W^{1,2}_{T_{(1)}},
 ~
 \Delta T_{(1,2)}\!:=\!T_{(1,2)}-T_{(1)},
 $$
so that
$$
 \begin{array}{rcl}
 u(0,x)
 &\!\!\!\!=&\!\!\!\!
 \E_{0,x}\Big[\1_{\{T_{(1)}=T\}}\frac{g(W_T)}{\bar F(T)}
 ~+~
 \1_{\{T_{(1)}<T\}} \frac{\Wc^{1}}{\rho(T_{(1)})}
 \\
 &&\hspace{15mm}
 \prod_{i=1}^2\!\!\Big( \1_{\{T_{(1,i)}=T\}}\frac{g(W^{1,i}_T)}{\bar F(\Delta T_{(1,i)})}
+  \1_{\{T_{(1,i)}<T\}}
\frac{(u D^{I^{1,i}}\!u)\big(T_{(1,i)},W^{1,i}_{T_{(1,i)}}\big)}
{\rho(\Delta T_{(1,i)})}
\Big)
              \Big].
 \end{array}
 $$
Our main representation result is obtained by iterating the last procedure, and solving the integrability problems which arise because of the singularity introduced by the random variable $\Wc^1$. The automatic differentiation, which is the main additional ingredient to the branching diffusion representation, is illustrated in the previous example when the operator $\Lc$ corresponds to the Brownian motion. This extends to the case of a more general diffusion operator by the so-called Bismuth-Elworthy-Li formula based on the Malliavin integration by parts formula, see Fourni\'e et al. \cite{FLLLT} for its use in the context of Monte Carlo approximation and the extension to other sensitivities.

Our main result provides a probabilistic representation of the solution of the semilinear PDE, with polynomial nonlinearity, in terms of a branching diffusion.
This requires naturally a technical condition ensuring the existence of a non-exploding solution for the PDE which can be either interpreted as a small maturity or a small nonlinearity condition. This new representation provides a new ingredient for the analysis of the corresponding PDE as it can be used to argue about existence, uniqueness, and regularity. We shall indeed prove a $C^1-$regularity result in order to prove the main Theorem \ref{thm:main1}. 

Moreover, our new representation has an important numerical implication as it is suitable for high dimensional Monte Carlo approximation. This is in fact the first high dimensional general method for nonlinear PDEs~! The practical performance of the method is illustrated on a numerical example in dimension $d=20$. The convergence of the numerical method is a direct consequence of the law of large numbers. The rate of convergence is also a direct consequence of the central limit theorem, and is therefore dimension-free. 
 The complexity of the method is easily shown to be of the order of $d^2$,
which cannot be avoided by the very nature of the equation whose second order term involves $d\times d$ matrices calculations.

The paper is organized as follows. Section \ref{sect:branchingdiff} introduces the marked branching diffusion. The main representation result is stated in Section \ref{sect:mainresult}. We next provide further discussions in Section \ref{sect:further} on the validity of our representation for systems of semilinear PDEs, and the possible combination with the unbiased simulation technique of \cite{HTT2,DOW}. The Monte Carlo numerical implications of our representation in high dimension are reported in Section \ref{sec:MC} with an illustration by a numerical experiment in dimension 20. Finally, we provide more numerical examples in Section \ref{sec:Simu2}.

\section{The marked branching diffusion}
\label{sect:branchingdiff}

\subsection{Semilinear PDE with polynomial nonlinearity}

	Let $d \ge 1$, $\M^d$ denotes the set of all $d \x d$ matrices,
	and $(\mu, \sigma): [0, T] \x \R^d \to \R^d \x \M^d$ the coefficient functions.
	For a function $u: [0,T]\x\R^d \to \R$, we denote by $Du$ and $D^2u$ the gradient and the Hessian of the function $u(t,x)$ w.r.t. variable $x$.
	Let $m \ge 0$ be a positive integer, we consider a subset $L \subset \N^{m+1}$,
	and a sequence of functions $(c_{\ell})_{\ell \in L}$ and $(b_i)_{i = 1, \cdots, m}$, 
	where $c_{\ell}: [0,T] \x \R^d \to \R$ and $b_i : [0,T] \x \R^d \to \R^d$.
	For every $\ell = (\ell_0, \ell_1, \cdots, \ell_m) \in L$, denote $| \ell | := \sum_{i=0}^m \ell_i$.
	A generator function $f: [0,T] \x \R^d \x \R \x \R^d$ is then defined by
	\be \label{eq:def_f}
		f(t,x,y, z) 
		&:=& \!\!
		\sum_{\ell = (\ell_0, \ell_1, \cdots, \ell_m) \in L} 
		c_{\ell}(t,x) ~y^{\ell_0} 
		~\prod_{i=1}^m
		\big(b_i(t,x) \cdot z \big)^{\ell_i}.
	\ee
	Given two matrix $A, B \in \M^d$, denote $A:B := \mbox{Trace}(AB^{\top})$.
	We will consider the following semilinear PDE:
	\be \label{eq:PDE}
		\partial_t u + \mu \cdot Du + \frac{1}{2} \sigma \sigma^{\top} \!:\! D^2 u + f(\cdot, u, Du)
		= 0,
		~\mbox{on}~[0,T) \x \R^d,
		~\mbox{and}~u(T,.)=g,
	\ee
	for some bounded Lipschitz function $g:\R^d\longrightarrow\R$.
	
	\begin{Remark}
	{\rm
		The nonlinearity \eqref{eq:def_f} includes the simplest case of a source term.
		Indeed, for $\ell = (0, 0, \cdots, 0)$, we have $c_{\ell}(t,x) y^{\ell_0} \prod_{i=1}^m \big(b_i(t,x) \cdot z \big)^{\ell_i}=c_{\ell}(t,x).$
	}
	\end{Remark}

\subsection{Age-dependent marked branching process}
\label{subsubsec:age_depend_proc}

	In preparation of the representation result, let us first introduce a branching process,
	characterized by a distribution density function $\rho: \R_+ \to \R_+$,
	a probability mass function $(p_{\ell})_{\ell \in L}$ (i.e. $p_{\ell} \ge 0$ and $\sum_{\ell \in L} p_{\ell} = 1$).

	Instead of the usual exponential arrival time, we shall consider a branching particle process with arrival time of distribution density function $\rho$. At the arrival time, the particle branches into $|\ell |$ offsprings with probability $p_{\ell}$,
	among which, $\ell_i$ particles carry the mark $i$, $i=0,\ldots,m$.
	Then regardless of its mark, each descendant particle performs the same but independent branching process as the initial particle.
	
	To construct the above process, 
	we will consider a probability space $(\Om, \Fc, \P)$ equipped with 
	\begin{itemize}
		\item a sequence of i.i.d. positive random variables $(\tau^{m,n})_{m, n \ge 1}$ of density function $\rho$, 
	
		\item a sequence of i.i.d. random elements $(I^{m,n})_{m, n \ge 1}$ with $\P(I^{m,n} = \ell) = p_{\ell}$, $\ell \in L$.
	\end{itemize}
	In addition, the sequences $(\tau^{m,n})_{m,n \ge 1}$ and $(I^{m,n})_{m,n \ge 1}$ are independent.
	

	We now construct an age-dependent branching process, with $(\tau^{m,n})_{m,n \ge 1}$ and $(I^{m,n})_{m,n \ge 1}$,
	using the following procedure.
	\begin{enumerate}
		\item We start from a particle marked by $0$, indexed by $(1)$, of generation $1$,
		whose arrival time is given by $T_{(1)} := \tau^{1,1} \wedge T$.
		
		\item Let $k = (k_1, \cdots, k_{n-1}, k_n) \in \N^n$ be a particle of generation $n$, with arrival time $T_k$.
		When $T_k < T$, we let $I_k =  I^{n, \pi_n(k)}$, where
		$$
			\pi_n ~\mbox{is an injection from}~ \N^n ~\mbox{to}~ \N,
		$$
		and at time $T_k$, it branches into $| I_k |$ offspring particles,
		which constitute $n+1$-the generation, and are indexed by $(k_1, \cdots, k_n, i)$ for $i = 1, \cdots , |I_k|$.

		\item \label{item:step_theta_k}
		When $I_k = (\hat \ell_0, \hat\ell_1, \cdots, \hat\ell_m)$, we have $| \hat \ell|$ offspring particles,
		among which we mark the first $\hat \ell_0$ particles by $0$,
		the next $\hat \ell_1$ particles by $1$, and so on,
		so that each particle has a mark $i$ for $i = 0, \cdots, m$.

		\item For a particle $k = (k_1, \cdots, k_n, k_{n+1})$ of generation $n+1$, 
		we denote by $k- := (k_1, \cdots, k_n)$ the ``parent'' particle of $k$,
		and the arrival time of $k$ is given by $T_k := \big(T_{k-} + \tau^{n+1, \pi_{n+1}(k)} \big) \wedge T$.
		
		\item In particular, for a particle $k = (k_1, \cdots, k_n)$ of generation $n$,
		and $T_{k-}$ is its birth time and also the arrival time of $k-$.
		Moreover, for the initial particle $k = (1)$, one has $k- = \emptyset$, and $T_{\emptyset} = 0$.
	\end{enumerate}
	
	The above procedure defines a marked age-dependent branching process.
	We denote further
	$$
		\theta_k := \mbox{mark of}~k,
		~~~
		\Kc^n_t := \begin{cases}
			\big\{ k ~\mbox{of generation}~n~\mbox{s.t.}~T_{k-} \le t < T_k \big\}, &\mbox{when}~~t \in [0,T),\\
			\{k ~\mbox{of generation}~ n~\mbox{s.t.}~ T_k = T\}, &\mbox{when}~~ t = T,
		\end{cases}
	$$
	and also
	$$
		\Kcb^n_t := \cup_{s \le t} \Kc^n_s,
		~~~~
		\Kc_t := \cup_{n \ge 1} \Kc^n_t
		~~~\mbox{and}~~~~
		\Kcb_t := \cup_{n \ge 1} \Kcb^n_t.
	$$
	Clearly, $\Kc_t$ (resp. $\Kc^n_t$) denotes the set of all living particles (resp. of generation $n$) in the system at time $t$,
	and $\Kcb_t$ (resp. $\Kcb^n_t$) denotes the set of all particles (resp. of generation $n$) which have been alive before time $t$.

	\begin{Example}
		Let us consider the case $d= 1$, with
		$$
			f(t,x,y,z) ~:=~ c_{0,0}(t,x) ~+~ c_{1,0}(t,x) y ~+~ c_{1,1}(t,x) y z .
		$$
		In this case,  $m = 1$, $L = \{\bar \ell_1 = (1,0), \bar \ell_2 = (1,1)\}$.
		For the sake of clarity, we present an typical path of the associated age-dependent process,
		with graphical illustration below. 
		The process starts from time $0$ with one particle indexed by $(1)$. 
		At terminal time $T$, the number of particles alive is $3$, with
		\b*
			\Kc_T 
			&=& 
			\big\{ (1,2,1), (1,1,1,1), (1,1,1,2) \big\},
		\e*
		\begin{itemize}
			\item At time $T_{(1)}$, particle $(1)$ branches into two particles $(1,1)$ and $(1,2)$.
			\item At time $T_{(1,1)}$, particle $(1,1)$ branches into $(1,1,1)$ and $(1,1,2)$.
			\item At time $T_{(1,2)}$, particle $(1,2)$ branches into $(1,2,1)$.
			\item At time $T_{(1,1,2)}$, particle $(1,1,2)$ dies out without any offspring particle.
			\item At time $T_{(1,1,1)}$, particle $(1,1,1)$ branches into $(1,1,1,1)$ and $(1,1,1,2)$.
			\item The particles in blue are marked by $0$, and the particles in red are marked by $1$.
		\end{itemize}
	\end{Example}

                \setlength{\unitlength}{0.8cm}
                \begin{picture}(20,10)
                        \thicklines

                        \put(0.5,0.5){\vector(1,0){15}}

                        \put(0.92,0){ $0$}
                        \put(1,0.5){\line(0,1){0.15}}

                        \put(3.72,0){ $T_{(1)}$}
                        \put(4,0.5){\line(0,1){0.15}}

                        \put(6.72,0){ $T_{(1,1)}$}
                        \put(7,0.5){\line(0,1){0.15}}

                        \put(8.22,0){ $T_{(1,2)}$}
                        \put(8.5, 0.5){\line(0,1){0.15}}

                        \put(10.5 ,0){ $T_{(1,1,2)}$}
                        \put(11, 0.5){\line(0,1){0.15}}

                        \put(11.85,0){ $T_{(1,1,1)}$}
                        \put(12, 0.5){\line(0,1){0.15}}

                        \put(13.92,0){ $T$}
                        \put(14, 0.5){\line(0,1){0.15}}

                        \put(1, 5){\cblue \line(1,0){3}}

                        \put(4, 5){\cblue \line(2,1){3}}
                        \put(4, 5){\cred \line(2,-1){4.5}}

                        \put(7, 6.5){\cblue \line(5,2){5}}
                        \put(7, 6.5){\cred \line(3,-1){4}}

                        \put(8.5, 2.75){\cblue \line(5,-1){5.5}}


                        \put(12, 8.5){\cblue \line(3,1){2}}
                        \put(12, 8.5){\cred \line(5,-2){2}}

                        \put(1.1, 4.6){{\scriptsize $(1)$}}

                        \put(4.2, 4.15){{\scriptsize $(1,2)$}}
                        \put(4.2, 5.6){{\scriptsize $(1,1)$}}

                        \put(7.2, 7.15){{\tiny $(1,1,1)$}}
                        \put(7.2, 5.7){{\tiny $(1,1,2)$}}

                        \put(8.7, 3.15){{\tiny $(1,2,1)$}}


                        \put(12.2, 9.35){{\tiny $(1,1,1,1)$}}
                        \put(12.2, 7.45){{\tiny $(1,1,1,2)$}}

		\end{picture}
		
		\vspace{2mm}

		\begin{Proposition}
			Assume that $\sum_{\ell\in L}|\ell|p_\ell < \infty$.
			Then the age-dependent branching process is well defined on $[0,T]$, i.e. 
			the number of particles in $\Kcb_t$ is finite a.s. for all $t \in [0,T]$.
		\end{Proposition}
		\proof See e.g. Theorem 1 of Athreya and Ney \cite[Chapter IV.1]{AthreyaNey},
		or Harris \cite[pp. 138-139]{Harris}.
		\qed

\subsection{The marked branching diffusion}

We next equip each particle with a Brownian motion in order to define a branching Brownian motion. 
	
We consider a sequence of independent $d$-dimensional Brownian motion $(W^{m,n})_{m, n \ge 1}$, which are also independent of $(\tau^{m,n},I^{m,n})_{m, n \ge 1}$. Define $W^{(1)}_t = \Delta W^{(1)}_t := W^{1,1}_t$ for all $t \in \big[0, T_{(1)} \big]$ and
	then for each $k = (k_1, \cdots, k_n) \in \Kcb_T \setminus \{(1)\}$,
	define 
	\be \label{eq:def_Wk}
		W^k_t ~:=~ W^{k-}_{T_{k-}} + \Delta W^k_{t - T_{k-}},
		~~\mbox{with}~~
		\Delta W^k_{t - T_{k-}} := W^{n, \pi_n(k)}_{t - T_{k-}}, ~~\mbox{for all}~ t \in [T_{k-}, T_k].
	\ee
	Then $(W^k_{\cdot})_{k \in \Kcb_T}$ is a branching Brownian motion.
	For each $k \in \Kcb_T$, we define an associated diffusion process $(X^k_t)_{t \in [T_{k-}, T_k]}$ 
	by means of the following SDE
	\be \label{eq:def_Xk}
		X^k_t ~=~ X^{k-}_{T_{k-}} + \int_{T_{k-}}^t \mu \big( s,  X^k_s \big) ds 
		+ \int_{T_{k-}}^t \sigma \big(s, X^k_s \big) dW^k_s,
		~~t \in [T_{k-}, T_k],~~\P\mbox{-a.s.,}
	\ee
	where for particle $(1)$, we fix the initial condition $X^{(1)}_0 = x_0$ for some constant $x_0 \in \R^d$. The well-posedness of the last SDE is guaranteed by standard conditions on the coefficients $\mu,\sigma$ contained in Assumption \ref{assum:GWtree}.
	
	The process $(X^k_{\cdot})_{k \in \Kcb_T}$ is our main marked branching diffusion process.	
	We finally introduce the sub-$\sigma$-fields
	\be \label{eq:Fc0}
		\Fc_0 := \sigma \big\{ \tau^{m,n}, I^{m,n} : m,n \ge 1 \big\},
		~
		\Fc_m := \sigma \big\{ W^{i,n}, \tau^{i,n}, I^{i,n} : n \ge 1, i \le m \big\},
		~m \ge 1.~~
	\ee

\section{The main representation} 
\label{sect:mainresult}
	We shall provide a representation result for a class the semilinear PDEs \eqref{eq:PDE} under general abstract conditions. More explicit sufficient conditions are provided later.

\subsection{Branching diffusion representation of semilinear PDEs}

	We first collect the conditions on the marked branching diffusion which are needed for our main results.

	\begin{Assumption} \label{assum:GWtree}
		\rmi The probability mass function $(p_{\ell})_{\ell\in L}$ satisfies 
		$p_{\ell} >0$ for all $\ell \in L$, and $\sum_{\ell \in L} |\ell |~ p_{\ell}  < \infty$.
		The density function $\rho: \R_+ \to \R_+$ is continuous and strictly positive on $[0, T]$, and such that
		$
			\Fb(T) :=\int_T^{\infty} \rho(t) dt > 0.
		$
		
		\noindent \rmii  $(\mu,\sigma): [0,T] \x \R^d \to \R^d \x \M^d$ are bounded continuous, and Lipschitz in $x$.
		
		\noindent \rmiii $c_{\ell}: [0,T] \x \R^d \to \R$ and $b_i : [0,T] \x \R^d \to \R^d$ are bounded continuous.
	\end{Assumption}

	Our next assumption is the key automatic differentiation condition on the underlying diffusion $\Xb^{t,x}_s$ defined by
	\be \label{eq:SDE_Xb}
		\Xb^{t,x}_s ~=~ x + \int_t^s \mu \big(r, \Xb^{t,x}_r \big) dr + \int_t^s \sigma \big(r, \Xb^{t,x}_r \big) dW_r,
		~~~s \in [t, T],
	\ee
	where $W$ is a $d$-dimensional Brownian motion independent of the branching diffusion.
	\begin{Assumption} \label{assum:auto_diff}
		There is a measurable functional $\Wcb(t, s, x, (W_r - W_t)_{r \in [t, s]})$ satisfying
		$(t,x) \mapsto \Wcb(t, s, x, (W_r - W_t)_{r \in [t, s]})$ is continuous,
		and for any $s \in [t, T]$ and bounded measurable function $\phi: \R^d \to \R$, one has
		$$
			\partial_x \E\big[ \phi \big( \Xb^{t,x}_s \big) \big]
			~=~
			\E \Big[ \phi\big( \Xb^{t,x}_s \big) \Wcb(t, s, x, (W_r - W_t)_{r \in [t, s]}) \Big].
		$$
	\end{Assumption}
	\begin{Remark}
		In case $(\mu, \sigma) \equiv (\mu_0, \sigma_0)$ for some constant $(\mu_0, \sigma_0) \in \R^d \x \M^d$, where $\sigma_0$ is not generate,
		then an example of such automatic differentiation function can be given by
		$$
			\Wcb \big(t, s, x, (W_r - W_t)_{r \in [t, s]} \big) ~:=~ (\sigma_0^{\top})^{-1} \frac{W_s - W_t}{s -t }.
		$$
		For general coefficient functions $(\mu, \sigma)$ satisfying some regularity and non-degeneracy conditions, one can find such functional $\Wcb$ using Malliavin calculus (see more discussions in Section \ref{subsubsec:ESC}).
	\end{Remark}

	Now, for each particle $k \in \Kcb_T$, we recall that it is born at time $T_{k-}$ and dies out at time $T_k$,
	its mark is given by $\theta_k$ and its branching type is given by $I_k$.
	Let us denote
	\be \label{eq:def_Wk}
		\Wc_k ~:=~ 
		\1_{\{\theta_k = 0\}}
		+
		\1_{\{\theta_k \neq 0\}}
		~b_{\theta_k}(T_{k-}, X^k_{T_{k-}})
		\cdot \Wcb \big( T_{k-}, T_k, X^k_{T_{k-}}, \Delta W^k_{\cdot}  \big).
	\ee
	We next introduce for a smooth function $u \in C^{1,2}([0,T]\x\R^d)$
	\be \label{eq:def_psi_n}
		\psi_n
		&:=&
		\Big[ \prod_{k \in \cup_{j=1}^n \Kc^j_T} \frac{g(X^k_T) - g(X^k_{T_{k-}}) \1_{\{\theta_k \neq 0\}} }{ \Fb(\Delta T_k)}  \Wc_k \Big]
		\Big[\!\!\! \prod_{k \in \cup_{j=1}^n (\Kcb^j_T \setminus \Kc^j_T)}
		\frac{c_{I_k}(T_k, X^k_{T_k})}{p_{I_k}} \frac{\Wc_k}{\rho(\Delta T_k) } \Big] \nonumber \\
		&&
		\Big[
			\prod_{k \in \Kcb^{n+1}_T} 
			\Big(\1_{\{\theta_k = 0\}} u +  \sum_{i=1}^m \1_{\{\theta_k = i \}} b_i \cdot Du \Big) (T_{k-}, X^k_{T_{k-}})
		\Big],
	\ee
	for all $n \ge 1$, and the corresponding limit
	\be \label{eq:def_psi}
		\psi 
		~:=~
		\Big[ \prod_{k \in \Kc_T} \frac{g(X^k_T) - g(X^k_{T_{k-}}) \1_{\{\theta_k \neq 0\}} }{ \Fb(\Delta T_k)}  \Wc_k \Big]
		\Big[ \prod_{k \in \Kcb_T \setminus \Kc_T}
		\frac{c_{I_k}(T_k, X^k_{T_k})}{p_{I_k}} \frac{\Wc_k}{\rho(\Delta T_k) } \Big].
	\ee
	Further, notice that the above branching diffusion process $(T_k, X^k_{\cdot})_{k \in \Kcb_T}$ and random variables $\psi$, $\psi_n$ are defined with initial condition $(0, x_0)$ on interval $[0,T]$.
	By exactly the same way, we can define the system with initial condition $(t,x)$ on interval $[t,T]$, let us denote them respectively by $(T^t_k, W^{t,k}_{\cdot}, X^{t,x,k}_{\cdot})_{k \in \Kcb^t_T}$, $\psi^{t,x}$, and $\psi^{t,x}_n$.
	
	We now provide a first result, under strong regularity conditions, which provides a better understanding of our representation. We emphasize that our main representation result in Theorem \ref{thm:main1} below will be established under more transparent conditions.

	\begin{Proposition} \label{prop:main1}
		Let Assumptions \ref{assum:GWtree} and \ref{assum:auto_diff} hold true. Suppose that the PDE \eqref{eq:PDE} has a solution $u \in C^{1,2}([0,T] \x \R^d)$ with $\E \Big[ \int_t^T \big| f(\cdot, u, Du)(s, \Xb^{t,x}_s) \big| ds \Big] <\infty$, for some $(t,x) \in [0,T]\x\R^d$. Assume further that $(\psi^{t,x}_n)_{n \ge 1}$ is uniformly integrable. Then
		$$
			\psi^{t,x} \in \L^1 
			~~~\mbox{and}~~
			u(t, x) = \E\big[ \psi^{t,x} \big].
		$$
	\end{Proposition}
	\proof \rmi It suffices to consider $(t,x)=(0,x_0)$. Since $g$ is bounded, it follows from the integrability condition on the process $f(.,u,Du)$ and the Feynma-Kac formula that
	\be \label{eq:u_1}
	u(0,x_0)
	&=&
	\E \Big[
		\frac{1}{\Fb(T)} g \big(\Xb^{0,x_0}_{T} \big) \Fb(T)
		+
		\int_0^T \frac{1}{\rho(s)} f \big(\cdot, u, Du\big) \big(s, \Xb^{0,x_0}_s\big) \rho(s) ds
		\Big]
	\nonumber\\
	&=&
	\E \Big[\frac{1}{\Fb(T_{(1)})} g \big(X^{(1)}_{T} \big) \1_{\{T_{(1)} = T\}} ~+ \nonumber \\
		\!\!\!&&\!\!
		~~~~~
		\frac{1}{\rho( T_{(1)})}
		\Big(
		\frac{c_{I_{(1)}}}{p_{I_{(1)}}} 
		u^{I_{(1),0}} \prod_{i=1}^m \big( b_i \cdot D u \big)^{I_{(1),i}} \Big)
		\big(T_{(1)}, X^{(1)}_{T_{(1)}}\big) \1_{\{T_{(1)} < T\}}
		\Big]
		\\
		&=&
		\E[\psi_1].\nonumber
	\ee
	\rmii
	Next, let $b_0 \in \R^d$ be a constant vector, and assume in addition that 
	the random variable $\psi_1 \big( b_0 \cdot \Wcb(0, T_{(1)},x_0, \Delta W^{(1)}_{\cdot}) \big)$ is integrable.
	Then under Assumptions \ref{assum:auto_diff}, 
		\be \label{eq:Du_1}
		b_0 \cdot D_x u(0, x_0)
		\!\!\!&=&\!\!\!
		\E \Big[ \psi_1 ~b_0 \cdot \Wcb \big(0, T_{(1)}, x_0, \Delta W^{(1)}_{\cdot} \big)  \Big] \nonumber \\
		\!\!\!&=&\!\!\!
		\E\Big[
			\Big( \psi_1- \frac{1}{\Fb(T_{(1)})} g(x_0) \1_{\{T_{(1)} = T\}} \Big)
			~b_0 \cdot \Wcb \big(0, T_{(1)}, x_0, \Delta W^{(1)}_{\cdot} \big)
		\Big],~~~~~~
	\ee
	where the first equality follows by Lemma A.3 of \cite{HTT2} 
	and the second equality follows from the fact that 
	$\E \big[ \Wcb \big(0, T, x, (W_s)_{s \in [0,T]} \big) \big] = 0$.

	\vspace{1mm}

	\noindent \rmiii For $k \in \Kcb^2_T$, change the initial condition from $(0, x_0)$ to 
	$(T_{k-}, X^k_{T_{k-}}) = (T_{(1)}, X^{(1)}_{T_{(1)}})$ in formula \eqref{eq:u_1} and \eqref{eq:Du_1}. Then, with $\Fc_1$ defined in \eqref{eq:Fc0},		\b*
		u(T_{k-}, X^k_{T_{k-}}) 
		\!\!&=&\!\!
		\E \Big[ 
		 \1_{\{k \in \Kc^2_T\}} 
		 \frac{g(X^k_T)}{\Fb(\Delta T_k)}
		+
		\1_{\{ k \in \Kcb^2_T \setminus \Kc^2_T\}}
		 \Psi_k
		~\Big| \Fc_1 \Big],
	\e*
	by the Markov property, and by Assumption \ref{assum:auto_diff},
	\b*
		Du(T_{k-}, X^k_{T_{k-}}) 
		\!\!&=&\!\!
		\E \Big[ 
		\Big(
		 \1_{\{k \in \Kc^2_T\}} 
		 \frac{g(X^k_T) - g(X^k_{T_{k-}})}{\Fb(\Delta T_k)}
		+
		\1_{\{ k \in \Kcb^2_T \setminus \Kc^2_T\}} \Psi_k
		\Big) \\
		&&~~~~~~~~~~~~~~~~~~~~~~~~~~~~~~
		\Wcb(T_{k-}, T_k, X^k_{T_{k-}}, \Delta W^k_{\cdot})
		~\Big| \Fc_1 \Big],
	\e*
	where $
		\Psi_k \!:=\! 
		\frac{1}{\rho(\Delta T_k)}  \frac{c_{I_k}(T_k, X^k_{T_k})}{p_{I_k}}
		\Big[ \prod_{k'-=k} 
			\Big(\1_{\{\theta_{k'} = 0\}} u +  \sum_{i=1}^m \1_{\{\theta_{k'} = i \}} b_i \cdot Du \Big) (T_{k'-}, X^{k'}_{T_{k'-}})
		\Big].
	$
	Plugging these expressions in the definition of $\psi_1$ in \eqref{eq:def_psi_n}, it follows from the integrability of $\psi_2$ and the tower property of conditional expectations that $u(0, x_0) = \E \big[ \psi_2 \big]$.

	\noindent \rmiv Iterating this procedure, we see that
	$$
		u(0, x_0) = \E \big[ \psi_n \big], ~~\mbox{for all}~ n \ge 1,
		~~~\mbox{and therefore}~~ u(0,x_0) = \lim_{n \to \infty} \E\big[ \psi_n \big] = \E \big[ \psi \big],
	$$
	where the last equality follows by the uniform integrability condition of $(\psi_n)_{n\ge1}$.
	\qed
	
\vspace{5mm}

	We now state our main representation result under abstract conditions on the automatic differentiation weight function $\Wcb(\cdot)$ involving the slight modification of $\psi$:
	\be \label{eq:def_psi_t}
		\tilde \psi 
		:=
		\Big[\!\! \prod_{k \in \Kc_T}\!\! \frac{g(X^k_T) - g(X^k_{T_{k-}}) \1_{\{\theta_k \neq 0 ~\mbox{or}~ k = (1) \}} }{ \Fb(\Delta T_k)}  \Wc_k \Big]
		\Big[\!\! \prod_{k \in \Kcb_T \setminus \Kc_T}
		\!\!\frac{c_{I_k}(T_k, X^k_{T_k})}{p_{I_k}} \frac{\Wc_k}{\rho(\Delta T_k) } \Big],
	\ee
with $\tilde\psi^{t,x}$ defined by an obvious change of origin. Explicit sufficient conditions for the validity of the next result will be reported in Section \ref{subsubsec:ESC} below. 
	
	\begin{Theorem} \label{thm:main1}
		Let Assumptions \ref{assum:GWtree} and \ref{assum:auto_diff} hold true, and suppose in addition that
		for all $(t,x) \in [0,T] \x \R^d$, there is some $\eps > 0$ such that
		$$
			(\psi^{s,y})_{(s,y) \in B_{\eps}(t,x)}
			~~~\mbox{and}~~
			\big(\tilde \psi^{s,y} \Wcb(s, T^s_{(1)}, y, \Delta W^{s,(1)}_{\cdot}) \big)_{(s,y) \in B_{\eps}(t,x)}
		$$ 
		are uniformly integrable,
		where $B_{\eps}(t,x):=\{(s,y) \in [0,T] \x \R^d : |s-t| + |x-y| \le \eps\}$. Then, the function $u(t,x) := \E[ \psi^{t,x}]$ is a continuous viscosity solution of the semilinear PDE \eqref{eq:PDE}. Moreover $Du$ exists and is continuous.
	\end{Theorem}
	\proof 
	\rmi Notice that the solution of SDE \eqref{eq:SDE_Xb} is continuous w.r.t. its initial condition $(t,x)$, and recall that $(t,x) \mapsto \Wcb(t,s,x, (W_r -W_t)_{r \in [t,s]}$ is also continuous,
	then under the uniform integrability condition on $(\psi^{t,x})$,
	one obtains that $u:[0,T] \x \R^d \to \R$ is continuous.
	Similarly, let us define 
	$$
		v_i(t,x) 
		:=
		\E \big[ \tilde \psi^{t,x} b_i(t,x) \cdot \Wcb(t, T^t_{(1)}, x, \Delta W^{t,(1)}_{\cdot}) \big]
		=
		\E \big[ \psi^{t,x} b_i(t,x) \cdot \Wcb(t, T^t_{(1)}, x, \Delta W^{t,(1)}_{\cdot}) \big],
	$$
	which is also continuous by the uniformly integrability condition.

	\vspace{1mm}

	\noindent \rmii Let us define $\phi :[0,T] \x \R^d \to \R$ by
	\be \label{eq:phi_expression}
		\phi \big(I_{(1)}, T_{(1)}, X^{(1)}_{T_{(1)}} \big) 
		\!\!\!&:=&\!\!\!
		\frac{1}{\Fb(T_{(1)})} g \big(X^{(1)}_{T} \big) \1_{\{T_{(1)} = T\}} \\
		\!\!\!\!\!\!\!\!&&\!\!\!\!\!
		+
		\frac{1}{\rho( T_{(1)})}
		\Big(
		\frac{c_{I_{(1)}}}{p_{I_{(1)}}} 
			\big(\psi^{T_{(1)}, X^{(1)}_{T_{(1)}}} \big)^{I_{(1),0}} 
			\prod_{i=1}^m v_i \big(T_{(1)}, X^{(1)}_{T_{(1)}} \big)^{I_{(1),i}} \Big)
		\1_{\{T_{(1)} < T\}} \nonumber \\
		&=&
		\E \big[ \psi ~\big|~ \Fc_1 \big], \nonumber
	\ee
	where $\Fc_1$ is defined in \eqref{eq:Fc0}.
	Notice that $\phi(i,t,x)$ is continuous in $x$, 
	then it follows by Assumption \ref{assum:auto_diff} and Lemma A.3 of \cite{HTT2} that
	$$
		D u(0,x_0) 
		= 
		\E\big[ \phi\big(I_{(1)}, T_{(1)}, X^{(1)}_{T_{(1)}} \big)  \Wcb(0, T_{(1)}, x_0, \Delta W^{(1)}_{\cdot}) \big]
		=
		\E \big[ \psi \Wcb(0, T_{(1)}, x_0, \Delta W^{(1)}_{\cdot}) \big].
	$$
	By changing the initial condition from $(0, x_0)$ to $(t,x)$
	and notice that 
	$$\E[\Wcb(t, T^t_{(1)}, x, \Delta W^{t,(1)}_{\cdot}) \1_{\{T^t_{(1)} = T\}}] = 0,$$
	it follows that
	$$
		Du(t,x) ~=~ \E \big[ \tilde \psi^{t,x} \Wcb(t, T^t_{(1)}, x, \Delta W^{t,(1)}_{\cdot}) \big],
	$$
	and one obtains that $Du: [0,T] \x \R^d \to \R^d$ is continuous from the uniform integrability of
	$\big(\tilde\psi^{t,x} \Wcb(t, T^t_{(1)}, x, \Delta W^{t,(1)}_{\cdot}) \big)$.
	Moreover, one has $v_i(t,x) = b_i(t,x) \cdot Du(t,x)$.
	
	\vspace{1mm}
	
	\noindent \rmiii Using the expression in \eqref{eq:phi_expression} and the law of $I_{(1)}$ and $T_{(1)}$, 
	and with similar arguments as in \eqref{eq:u_1},
	it follows that
	\b*
		u(t,x) ~=~ \E[ \psi^{t,x}] 
		&=& 
		\E \Big[ g\big(\Xb^{t,x}_T \big) + \int_t^T f \big(\cdot, u, Du \big) (s, \Xb^{t,x}_s)  ds \Big].
	\e*
	Let $h > 0$, denote $\mbox{\sc h}_h  :=  (t+h) \wedge \inf\{s >t ~: |\Xb^{t,x}_s -x | \ge 1\}$, 
	then by the flow property of $\Xb^{t,x}$, one has
	$$
		u(t,x) ~=~
		\E\Big[ u(\mbox{\sc h}_h, \Xb^{t,x}_{\mbox{\sc h}_h}) + \int_t^{\mbox{\sc h}_h} f \big(\cdot, u, Du \big) (s, \Xb^{t,x}_s)  ds
		\Big],
	$$
	and we may verify by standard arguments that $u$ is a viscosity solution of PDE \eqref{eq:PDE}.
	\qed

\subsection{More explicit sufficient conditions}
\label{subsubsec:ESC}

	We now provide some explicit sufficient conditions  which guarantee the validity of the conditions of Theorem \ref{thm:main1}.
	Define $| \varphi |_{\infty} := \sup_{x \in \R^d} |\varphi (x)|$ for any bounded function $\varphi: \R^d \to \R$,
	and $|\phi|_{\infty} := \sup_{i = 1}^d |\phi_i|_{\infty}$ for any bounded vector function $\phi: \R^d \to \R^d$.

	We first recall the Bismut-Elworthy-Li formula from Malliavin calculus,
	which was used by Fourni\'e, Lasry, Lebuchoux, Lions and Touzi \cite{FLLLT} as an automatic differentiation tool, 
	see also \cite{BET}, \cite{BouchardTouzi} and \cite{FTW} for subsequent usefulness of the automatic differentiation in the context of the Monte Carlo approximation of nonlinear PDEs.
	We emphasize that such automatic differentiation function is not unique.

	\begin{Assumption} \label{assum:diffussion_coef}
		The coefficients $\mu,\sigma$ are bounded continuous, with bounded continuous partial gradients $D \mu,D \sigma$, and $\sigma$ is uniformly elliptic.
	\end{Assumption}
	Notice that $(\Xb^{t,x}_s)_{s \in [t,T]}$, as defined by \eqref{eq:SDE_Xb},
	is completely determined by $(t, x, (W_s - W_t)_{s \in [t, T]})$.
	We then introduce the corresponding first variation process $Y$:
	\be \label{eq:def_Y}
		Y_t := \I_d,
		~
		d Y_s = D \mu(s, \Xb^{t,x}_s) Y_s ds + \!\sum_{i=1}^d D \sigma_i(s, \Xb^{t,x}_s) Y_s dW^i_s,
		~\mbox{for}~ s \in [t, T], \!~\P\mbox{-a.s.},~~
	\ee
	where $\I_d$ denotes the $d \x d$ identity matrix,
	and $\sigma_i(t,x) \in \R^d$ denotes the $i$-th column of matrix $\sigma(t,x)$.
	Then one has the following result (see e.g. Exercise 2.3.5 of Nualart \cite[p.p. 125]{Nualart}, or Proposition 3.2. of \cite{FLLLT}).
	\begin{Proposition} \label{prop:MalliavinWeight}
		Let Assumption \ref{assum:diffussion_coef} hold true, then Assumption \ref{assum:auto_diff} holds true 
		with the choice of automatic differentiation function $\Wcb$ defined by
		\be \label{eq:MalliavinWeight}
			\Wcb \big(t, s, x, (W_r - W_t)_{r \in [t, s]} \big)
			~:=~ 
			\frac{1}{s-t}\int_t^s \big[ \sigma^{-1}(r, \Xb^{t,x}_r) Y_r \big]^{\intercal} dW_r.
		\ee
	\end{Proposition}

	\begin{Remark} \label{rem:MWeight}
		When $\mu \equiv 0$ and $\sigma(t,x) \equiv \sigma_0$ for some non-degenerate constant matrix $\sigma_0 \in \M^d$, 
		one then has $Y_t \equiv \I_d$ and so that
		$$
			\Wcb \big(t, s, x, (W_r - W_t)_{r \in [t, s]} \big)
			~=~
			\big(\sigma_0^{\top} \big)^{-1}  \frac{W_s -W_t }{s - t}.
		$$
	
	\end{Remark}

	With the above choice of automatic differentiation weight function \eqref{eq:MalliavinWeight}, 
	we can now derive some upper bounds for random variables $(\psi_n, n \ge 1)$. 
	Recall that $L_g$ is the Lipschitz constant of $g$, denote by $B^{\infty}_0(L_g):= \{(x_1, \cdots, x_d) \in \R^d~: |x_i| \le L_g\}$
	and $\Wcb_{t,x,s} := \Wcb \big(t, s, x, (W_r - W_t)_{r \in [t, s]} \big)$.
	Then for $n \ge 1$, $q > 1$, we introduce two constants $C_{1,q}$ and $C_{2,q}$ by
	$$
		C_{1,q} := |g|_{\infty}^q 
		\vee 
		\sup_{0 \le t < s \le T, ~x \in \R^d, ~i = 1, \cdots, m, ~b_0 \in B^{\infty}_0(L_g)} 
		\E \Big[ \Big| \big( b_0 \cdot (\Xb^{t,x}_s -x )\big) 
			\big( b_i(t,x) \cdot \Wcb_{t,x,s} \big) \Big|^q 
		\Big]
	$$
	and
	$$
	 	C_{2,q} 
		:=
		\sup_{0 \le t < s \le T, ~x \in \R^d, ~i = 1, \cdots, m}
		\E \Big[ \big| \sqrt{s-t} ~b_i(t,x) \cdot \Wcb_{t,x,s} \big|^q \Big],
	$$
	and then
	$$
		\widehat C_{1,q} ~:=~ \frac{C_{1,q} }{ \Fb(T)^{q-1}},
		~~~~~~
		\widehat C_{2,q} ~:=~ C_{2,q}
			\sup_{\ell \in L, ~t \in (0, T]}
			\Big(\frac{|c_{\ell}|_{\infty}}{p_{\ell}} \frac{ t^{- \frac{q}{2(q-1)}}}{\rho(t)} \Big)^{q-1}.
	$$

	\begin{Remark}
		\rmi Under Assumption \ref{assum:diffussion_coef}, the tangent process $Y$ is defined by a linear SDE,
		which has finite moment of any order $q \ge 1$.
		Then the two constant $C_{1,q}$ and $C_{2,q}$ are both finite.
		And for all $k \in \Kcb_T$, one has
		\be \label{eq:def_Cq}
			\max\Big\{ |g|_{\infty}^q , ~\E \Big[ \big| ( D_g \cdot \Delta X_k) \Wc_k \big|^q ~\Big| \Fc_0 \Big] \Big\} \le C_{1,q},
			~~
			\E \Big[ \Big(\sqrt{\Delta T_k} | \Wc_k |\Big)^q \Big| \Fc_0 \Big] \le C_{2,q},~~
		\ee
		where the sub-$\sigma$-field $\Fc_0$ is defined in \eqref{eq:Fc0}.

		\noindent \rmii Notice that for a random variable $N \sim N(0,1)$ and non-negative integer $q \ge 0$, one has $\E[|N|^q] = 2^{\frac{q}{2}} \Gamma\big( \frac{q+1}{2} \big)/ \sqrt{\pi}$.
		Then if $(\mu, \sigma) \equiv (0, \sigma_0)$, for some constant $(\mu_0, \sigma_0) \in \R^d \x \M^d$, and $\Wcb$ as in Remark \ref{rem:MWeight},
		it follows by direct computation that
		$$
			C_{1,q}
			\le
			|g|_{\infty}^q  \vee
			\Big(
			\sup_{b_0 \in B^{\infty}_0(L_g)} \!\!\!
			\big( b_0^{\top} \sigma_0 \sigma_0^{\top} b_0 \big)
			+
			\max_{i = 1, \cdots, m}\!\!
			\| b_i^{\top} (\sigma_0 \sigma_0^{\top})^{-1} b_i\|_{\infty} \Big)^q
			2^{q - 1} \Gamma\Big( \frac{2q+1}{2} \Big)/ \sqrt{\pi},
		$$
		and
		$$
			C_{2,q} 
			=
			\max_{i = 1, \cdots, m} \| b_i^{\top} \big(\sigma_0 \sigma_0^{\top} \big)^{-1} b_i \|_{\infty}^{\frac{q}{2}}
			~2^{\frac{q}{2}} \Gamma\Big( \frac{q+1}{2} \Big)/ \sqrt{\pi}.
		$$
	\end{Remark}
	We are now ready for the main explicit sufficient conditions for the validity of the representation Theorem \ref{thm:main1}. Notice that the following conditions can be interpreted either as a small maturity or small nonlinearity restriction.
	
	\begin{Assumption} \label{assum:integrability_explicit}
		For some $q > 1$, one of the following two items holds true.

		\noindent \rmi Both $C_{1,q} \Big( \frac{1}{ \Fb(T)} \Big)^q$ and
		$\sup_{\ell \in L, t \in (0,T]} C_{2,q} \Big(\frac{|c_{\ell}|_{\infty}}{p_{\ell}} \frac{1}{\sqrt{t} \rho(t)} \Big)^q$ 
		are bounded by $1$.
		
		\noindent \rmii $T<\int_{\widehat C_{1,q}}^{\infty} \big(\widehat C_{2,q} \sum_{\ell \in L} |c_{\ell}|_{\infty} ~x^{|\ell|}\big)^{-1}dx$.
	\end{Assumption}

	\begin{Remark} \label{rem:integrability_explicit}
		\rmi To ensure that 
		$\sup_{\ell \in L, t \in (0,T]} C_{2,q} \Big(\frac{|c_{\ell}|_{\infty}}{p_{\ell}} \frac{1}{\sqrt{t} \rho(t)} \Big)^q$
		is bounded by $1$,
		it is necessary to choose $(p_{\ell})_{\ell \in L}$ such that $\frac{| c_{\ell}|_{\infty}}{p_{\ell}}$ is uniformly bounded,
		and to choose a density function such that $\rho(t) \ge C t^{-1/2}$.
		
		\noindent \rmii To ensure that $\widehat C_{2,q}$ is finite, one needs to choose the density function $\rho$ such that $\rho(t) \ge C t^{- \frac{q}{2(q-1)}}$,
		and hence it is necessary that $q \in (2, \infty)$ so that $\frac{q}{2(q-1)} \in (\frac{1}{2}, 1)$.
	\end{Remark}

	\begin{Theorem} \label{thm:main2}
		Consider the automatic differentiation function \eqref{eq:MalliavinWeight},
		and suppose that Assumptions \ref{assum:GWtree}, \ref{assum:diffussion_coef} and \ref{assum:integrability_explicit} hold true.
		
		\noindent \rmi Then Assumptions \ref{assum:auto_diff} holds, and 
		$\big(\psi^{t,x}, \tilde \psi^{t,x} \Wcb(t, T^t_{(1)}, x, \Delta W^t_{(1)}) \big)_{(t,x) \in [0,T] \x \R^d}$ is uniformly integrable. Consequently, $u(t,x) := \E[\psi^{t,x}]$ is a viscosity solution of PDE \eqref{eq:PDE}.
		
		\noindent \rmii If Assumption \ref{assum:integrability_explicit} holds with some $q \ge 2$, then $\E\big[|\psi^{t,x}|^2\big]<\infty$.
	\end{Theorem}
	\proof 
	\rmi First, using Proposition \ref{prop:MalliavinWeight}, it is clear that Assumption \ref{assum:auto_diff} holds true with the choice of automatic differentiation function in \eqref{eq:MalliavinWeight}.

	\vspace{1mm}
	
	\noindent \rmii Next, for $q \ge 1$, let us introduce
	\b*
		\chi^q_{\infty}
		\!\!&:=&\!\!\!
		\Big[\!\!\! \prod_{k \in \Kc_T}\!\!\! 
		C_{1,q} \Big( \frac{1}{ \Fb(\Delta T_k)} \Big)^q \Big]
		\Big[\!\!\! \prod_{k \in \Kcb_T \setminus \Kc_T}\!\!\!\!\!\!
		C_{2,q} \Big(\frac{|c_{I_k}|_{\infty}}{p_{I_k}} \frac{1}{\sqrt{\Delta T_k} \rho(\Delta T_k)} \Big)^q \Big].
	\e*
	By conditioning on $\Fc_0$, it follows from \eqref{eq:def_Cq}, together with direct computation, that
	\be \label{eq:psi_le_chi}
		\E[ |\psi|^q ] ~\le~ \E [\chi^q_{\infty}]
		~~~\mbox{and}~~
		\E \big[  \big|  \tilde \psi \Wcb(t, T^t_{(1)}, x, \Delta W^t_{(1)}) \big|^q \big]
		\le C \E [\chi^q_{\infty}],
	\ee
	for some constant depending only on the Lipschitz constant $L_g$.
	
	\noindent \rmiii When Assumption \ref{assum:integrability_explicit} \rmi holds true for some $q > 1$,
	then it is clear that $\E [ |\psi|^q ] \le 1$.
	Notice that the above argument is independent of the initial condition $(0, x_0)$,
	it follows that $\big(\psi^{t,x}, \tilde\psi^{t,x} \Wcb(t, T^t_{(1)}, x, \Delta W^t_{(1)}) \big)_{(t,x) \in [0,T] \x \R^d}$ is uniformly integrable.
	
	\vspace{1mm}
	
	\noindent \rmiv When Assumption \ref{assum:integrability_explicit} \rmii holds true for some $q > 1$.
	Consider the ODE on $[0,T]$:
	\b*
		\eta(T) ~=~ \widehat C_{1,q},
		&&
		\eta'(t) ~+~ \sum_{\ell \in L} \widehat C_{2,q}~ \|c_{\ell}\|_{\infty} ~\eta(t)^{|\ell|} ~=~ 0.
	\e*
	Under Assumption \ref{assum:integrability_explicit} (ii), it is clear that the above ODE admits a unique finite solution on $[0,T]$.
	We next introduce a sequence of random variables
	\b*
		\widehat \chi_n^q
		\!\!&:=&\!\!
		\Big[ \prod_{k \in \cup_{j=1}^n \Kc^j_T} \frac{\widehat C_{1,q} }{\Fb(\Delta T_k)}  \Big]
		\Big[ \prod_{k \in \cup_{j=1}^n (\Kcb^j_T \setminus \Kc^j_T)}
		\widehat C_{2,q} \frac{|c_{I_k}|_{\infty}}{p_{I_k}} \frac{1}{\rho(\Delta T_k)} \Big]
		\Big[ \prod_{k \in \Kcb^{n+1}_T}  \eta(T_{k-}) \Big],
	\e*
	and
	\b*
		\widehat \chi^q_{\infty}
		~:=~
		\lim_{n\to\infty} \widehat \chi_n^q
		~=~
		\Big[ \prod_{k \in \Kc_T} \frac{\widehat C_{1,q} }{\Fb(\Delta T_k)}  \Big]
		\Big[ \prod_{k \in \Kcb_T \setminus \Kc_T}
		\widehat C_{2,q} \frac{|c_{I_k}|_{\infty}}{p_{I_k}} \frac{1}{\rho(\Delta T_k)} \Big].
	\e*
	Then by the same arguments as in the proof of Proposition \ref{prop:main1}, it is easy to check that
	\b*
		\eta(0) 
		~=~ 
		\eta(T) + \int_0^T  \sum_{\ell \in L} \widehat C_{2,q}~ |c_{\ell}|_{\infty} ~\eta(t)^{|\ell|} dt
		~=~
		\E\big[ \widehat \chi^q_1 \big]
		~=~
		\E \big[ \widehat \chi^q_n \big], ~~~\forall n \ge 1;
	\e*
	and hence by direct computation, it follows that
	$$
		\E \big[ \big|\psi\big|^q \big] 
		~\le~ 
		\E \big[ \chi^q_{\infty} \big] 
		~\le~ 
		\E \big[ \widehat \chi^q_{\infty} \big] 
		~\le~ 
		\liminf_{n \to \infty} \E \big[ \widehat \chi^q_n \big] 
		~=~
		\eta(0) ~<~ \infty.
	$$
	Changing the origin from $(0, x_0)$ to $(t,x)$, we see that
	$$
		\sup_{(t,x) \in [0,T] \x \R^d} \E \big[ |\psi^{t,x}|^q \big] 
		~\le~
		\sup_{t \in [0,T]} \eta(t)
		~<~
		\infty, 
	$$
	and hence $\big(\psi^{t,x} \big)_{(t,x) \in [0,T] \x \R^d}$ is uniformly integrable.
	The same arguments using \eqref{eq:psi_le_chi} show that
	$\big(\tilde\psi^{t,x} \Wcb(t, T^t_{(1)}, x, \Delta W^t_{(1)}) \big)_{(t,x) \in [0,T] \x \R^d}$ is uniformly integrable.
\qed


\section{Further discussions}
\label{sect:further}

\paragraph{Representation of the PDE system}

	Let us consider a PDE system $(v_j)_{j=1, \cdots, n}$, 
	where for each $j$, $v_j: [0,T] \x \R^d \to \R$ satisfies
	\b*
		\partial_t v_j + \mu_j \cdot D v_j + \frac{1}{2} a_j : D^2 v_j + f(\cdot, v_1, \cdots, v_n, D v_1, \cdots, D v_n) ~=~ 0,
	\e*
	for some diffusion coefficient function $(\mu_j, a_j): [0,T] \x \R^d \longrightarrow \R^d \x \S^d$,
	and some polynomial function $f: [0,T] \x \R^d \x \R^n \x (\R^d)^n \longrightarrow \R$.
	Our methodology immediately applies to this context, and provides a stochastic representation for the solution of the above PDE system, by means of a regime-changed branching diffusions:
at every branching time, the independent offspring particles perform subsequently different branching diffusion regime.

\paragraph{Representation in view of unbiased simulation}
        \label{para:unbiased}
	With the same idea of proof, 
	we can also obtain an alternative representation result, with a frozen coefficient SDE in place of SDE \eqref{eq:def_Xk}.
	When the coefficient function $a \equiv a_0$ for some constant $a_0 \in \S^d$,
	this has significant  application in terms of Monte Carlo approximation,
	as it leads to a representation random variable which can be simulated exactly,
	while the branching diffusion process $X^k_{\cdot}$ in \eqref{eq:def_Xk} needs a time discretization technique and hence creates some discretization error in the simulation.
	Let us present this alternative representation formula in the case of constant diffusion coefficient case, 
	i.e. $a \equiv a_0 = \sigma_0 \sigma_0^{\intercal}$ for some non-degenerate constant matrix $\sigma_0 \in \M^d$.	

	Let $\Lh := L \cup \{ \partial \}$, where $\partial$ represents an artificial index;
	$\ph = (\ph_{\ell})_{\ell \in \Lb}$ be a  probability mass function and $(\Ih^{m,n})_{m,n \ge 1}$ be a sequence of i.i.d. random variables of distribution $\ph$, and independent of the sequences of i.i.d Brownian motion $(\Wh^{m,n})_{m,n \ge 1}$ and i.i.d positive random variable $(\Th^{m,n})_{m,n \ge 1}$ of density function $\rho$.
	Then following exactly the same procedure in Section \ref{subsubsec:age_depend_proc}, we can construct another age-dependent branching process, denoted by $(\Th_k)_{k \in \Kcbh_T}$ with branching type $\Ih_k := \Ih^{n, \pi_n(k)}$.
	Here, when $\Ih_k = (\hat \ell_0, \cdots, \hat \ell_m) \in L$, it produces $|\hat \ell|$ offspring particles, marked by $i = 0, \cdots, m$ exactly as in Step \ref{item:step_theta_k} in the construction of age-dependent process $\Kcb_T$ in Section \ref{subsubsec:age_depend_proc};
	when $I_k = \partial$, it produces only one offspring particle, marked by $m+1$.
	Then for every $k \in \Kcbh_T$, we equipped it with an independent Brownian motion $\Wh^{k}_{\cdot}$ as in \eqref{eq:def_Wk}.
	Next, let us define $\Xh^{(1)}_0 = x_0$, and subsequently for every $k \in \Kcbh_T$,
	\be \label{eq:def_Xh}
		\Xh^k_{\Th_k} ~:=~ \Xh^k_{\Th_{k-}} +~ \mu(\Th_{k-}, \Xh^k_{\Th_{k-}}) \Delta \Th_k ~+~ \sigma_0 \Delta \Wh^k_{\Delta \Th_k},
		~~~
		\mbox{with}~~
		\Xh^k_{\Th_{k-}} := \Xh^{k-}_{\Th_{k-}}.
	\ee
	For this case, the automatic differentiation functions take a particularly simple formula,
	which is compatible with the purpose of the unbiased simulation algorithm.
	Let us introduce
	\be \label{eq:def_Wch_k}
		\Wch_k 
		&:=& 
		\1_{\{\theta_k = 0\}} 
		~+~ 
		 b_{\theta_k}(\Th_{k-}, \Xh^k_{\Th_{k-}}) \cdot (\sigma_0^{\top})^{-1} \frac{\Delta \Wh^k_{\Delta \Th_k}}{\Delta \Th_k} \1_{\{ \theta_k \in \{ 1, \cdots, m\} \}} \\
		&& +~
		\Big( \mu(\Th_{k-}, \Xh^k_{\Th_{k-}}) - \mu(\Th_{(k-)-}, \Xh^{k-}_{\Th_{(k-)-}})\Big)
		\cdot (\sigma_0^{\top})^{-1} \frac{\Delta \Wh^k_{\Delta \Th_k}}{\Delta \Th_k}  \1_{\{\theta_k = \partial\}}. \nonumber
	\ee
	Finally, setting $c_{\partial} \equiv 1$,
	and replacing $(X^k_{\cdot}, \Wc_k, p_{I_k} )$ in the definition of $\psi$ and $\tilde \psi$ (in and below \eqref{eq:def_psi})
	by $(\Xh^k_{\cdot}, \Wch_k, \ph_{\Ih_k})$, we obtain
	\be \label{eq:def_psi_exact}
		\psih
		~:=~
		\Big[ \prod_{k \in \Kch_T} \frac{g(\Xh^k_T) - g(\Xh^k_{T_{k-}}) \1_{\theta_k \neq 0} }{ \Fb(\Delta T_k)}  \Wch_k \Big]
		\Big[ \prod_{k \in \Kcbh_T \setminus \Kch_T}
		\frac{c_{\Ih_k}(\Th_k, \Xh^k_{\Th_k})}{\ph_{\Ih_k}} \frac{\Wch_k}{\rho(\Delta \Th_k) } \Big],
	\ee
	and similarly $\tilde \psih$.

	Next, given a constant vector $\mu_0 \in \R^d$, we keep the same branching Brownian motion $(\widehat W^k_{\cdot} )_{k \in \widehat \Kcb_T}$, and then
	introduce another diffusion process $\Xh^{\mu_0, k}_{\cdot}$ by
	$$
		\Xh^{\mu_0, k}_{\Th_{(1)}} ~:=~ x_0 ~+~ \mu_0 \Delta \Th_{(1)} + \sigma_0 \Delta \Wh^{(1)}_{\Delta T_{(1)}},
	$$
	and the subsequent process $\Xh^{\mu_0,k}_{\Th_k}$ for $k \in \Kcbh_T \setminus \{(1)\}$ by the same induction relation as in \eqref{eq:def_Xh}.
	We then introduce $\Wch_k^{\mu_0}$ as in \eqref{eq:def_Wch_k} by replacing $\Xh^k$ by $\Xh^{\mu_0,k}$,
	and replacing $\mu(\Th_{(k-)-}, \Xh^{k-}_{\Th_{(k-)-}})$ by $\mu_0$ when $k=(1)$.
	Replacing $(\Xh^k, \Wch_k)$ by $(\Xh^{\mu_0,k}, \Wch^{\mu_0}_k)$ in \eqref{eq:def_psi_exact},
	it defines a new random variable $\psih^{\mu_0}$.
	Finally, by changing the initial condition $(0, x_0)$ and time interval $[0,T]$ to 
	$(t,x)$ and $[t,T]$, one obtains $\Wh^{t,k}$, $\Th^t_k$, $\psih^{t,x}$, $\tilde \psih^{t,x}$, $\psih^{\mu_0,t,x}$ etc.
	
	\begin{Proposition}
		Suppose that Assumptions \ref{assum:GWtree} holds true,
		and the semilinear PDE \eqref{eq:PDE} has uniqueness for bounded viscosity solution.
		Suppose in addition that for every $(t,x) \in [0,T] \x \R^d$, and $\mu_0$ lies in a neighborhood of $\mu(t,x)$, one has
		$$
			\psih^{\mu_0, t,x}  
			~~\mbox{and}~~
			\psih^{\mu_0, t,x} \Delta \Wh^{t,(1)}_{\Delta \Th^t_{(1)}}/\Delta \Th_{t,(1)}
			~~~\mbox{is integrable,}
		$$
		and the family of random variables 
		$$
			(\psih^{t,x})_{(t,x) \in [0,T] \x \R^d}
			~~~\mbox{and}~~
			\big(\tilde \psih^{t,x} \Delta \Wh^{t,(1)}_{\Delta \Th^t_{(1)}}/\Delta \Th_{t,(1)} \big)_{(t,x) \in [0,T] \x \R^d}
		$$ 
		are uniformly integrable with uniformly bounded expectation, 
		define $\hat u(t,x) := \E[ \psi^{t,x}].$
		Then the derivative $D \hat u$ exists, $\hat u$ and $D \hat u$ are both continuous; 
		and moreover, $u$ is the unique bounded viscosity solution of semilinear PDE \eqref{eq:PDE}.
	\end{Proposition}
	\noindent {\bf Sketch of proof.} 
	\rmi First, by the uniform integrability condition, $\hat u$ is bounded continuous.
	Let us introduce 
	$$ 
		\tilde u(\mu_0, t,x) := \E [\psih^{\mu_0, t,x} ]
		~~~\mbox{and}~~
		\hat v(t,x) 
		~:=~ 
		\E \Big[ \tilde \psih^{t,x} \Delta \Wh^{t,(1)}_{\Delta \Th^t_{(1)}}/\Delta \Th_{t,(1)} \Big].
	$$
	Notice that $\hat v$ is uniformly bounded and continuous.
	Recall that $W$ is a standard $d$-dimensional Brownian motion independent of the branching diffusion process,
	we also introduce
	$$
		\Xh^{t,x}_s ~:=~ x ~+~ \mu_0 (s-t) ~+~ \sigma_0 (W_s - W_t),~~s \in [t, T],
	$$
	where $\mu_0 \in \R^d$ is a constant vector in a neighborhood of $\mu(t,x)$.
	Then one obtains as in \eqref{eq:u_1} that
	\b*
		\tilde u(\mu_0, t,x) 
		\!\!\!&=&\!\!\!
		\E \Big[ \frac{1}{\Fb( \Th^t_{(1)})} g \big( \Xh^{t,x}_{\Th^t_{(1)}} \big) \1_{\{\Th^t_{(1)} = T\}} 
		+
		\frac{\1_{\{\Ih_{(1)} = \partial \}}}{\rho(\Th^t_{(1)}) \ph_{\partial}} \big( (\mu - \mu_0)\cdot \hat v \big) \big(\Th^t_{(1)}, \Xh^{t,x}_{\Th^t_{(1)}} \big) \\
		&&~~~~~~~~~+~~
		\frac{\1_{\{\Ih_{(1)} \neq \partial  \}}}{\rho(\Th^t_{(1)})} 
		\Big( \frac{c_{\Ih_{(1)}}}{\ph_{\Ih_{(1)}}}  \hat u^{\Ih_{(1),0}}  \prod_{i=1}^m (b \cdot \hat v)^{\Ih_{(1),i}} \Big)
			\big(\Th^t_{(1)}, \Xh^{t,x}_{\Th^t_{(1)}} \big) \1_{\{\Th^t_{(1)} < T \}}~\Big]\\
		\!\!\!&=&\!\!\!
		\E \Big[ g \big(\Xh^{t,x}_T \big) + \int_t^T \Big((\mu- \mu_0) \cdot \hat v + f( \cdot, \hat u, \hat v) \Big) 
			\big(s, \Xh^{t,x}_s \big) ds \Big].
	\e*
	By standard argument, $(t,x) \mapsto \tilde u(\mu_0, t,x)$ is a viscosity solution of
	$$
		-~ \partial_t u ~+~ \mu_0 \cdot Du ~+~ \frac{1}{2} a_0 : D^2 u 
		~+~ 
		(\mu - \mu_0) \cdot \hat v
		+
		\sum_{\ell \in L} c_{\ell} \hat u^{\ell_0} \prod_{i=1}^m (b_i \cdot \hat v)^{\ell_i}
		~=~
		0,
	$$
	with terminal condition $g$.
	Since $\hat u$ and $\hat v$ are bounded continuous, the above PDE has uniqueness for bounded viscosity solution, which induces that $\tilde u(\mu_0, t,x)$ is independent of $\mu_0$
	and $\hat u(t,x) = \tilde u(\mu_0, t,x)$ for $\mu_0$ in a neighborhood of $\mu(t,x)$.
	
	\vspace{1mm}
	
	\noindent \rmii We can then compute the derivative
	$D_x \tilde u(\mu_0, t,x)$ and then set $\mu_0 := \mu(t,x)$, 
	it follows that 
	$$
		D \hat u(t,x) ~=~ D_x \tilde u(\mu(t,x), t,x) ~=~ \hat v(t,x),
	$$
	which is also bounded continuous.
	This implies that $\hat u(t,x)$ is a viscosity solution of \eqref{eq:PDE},
	and we hence conclude the proof by uniqueness of the viscosity solution of \eqref{eq:PDE}.
	\qed
	
	\vspace{2mm}
	
	The integrability and square integrability of $\psih$ 
	can be analyzed in exactly the same way as in Theorem \ref{thm:main2}.
	We just notice that the above defined random variable $\psih$ 
	can be simulated exactly from a sequence of Gaussian random variable, discrete distributed random variables $\Ih^{m,n}$ and r.v. $\Th^{m,n}$ of distribution density function $\rho$.
	It is then in particular interesting to serve as a Monte-Carlo estimator for $u(0, x_0)$.


	
	
\paragraph{On the representation of fully nonlinear PDEs}

	Formally, one can also obtain a representation result for fully nonlinear PDE,
	using the same automatic differentiation functions of order 2.
	However, this raises a serious integrability problem which can not be solved by conditions as in Assumption \ref{assum:integrability_explicit}.
	To illustrate the main difficulty, let us consider the following PDE in the one-dimensional case $d = 1$:
	\be \label{eq:PDE_1}
		u(T,x) = g(x),
		~~~
		\partial_t u ~+~ \frac{1}{2} D^2 u ~+~ f_0(D^2u) ~=~ 0,
		~~\mbox{on}~[0,T] \x \R^d,
	\ee
	where $f_0(\gamma)  = c_0 \gamma$ for some constant $c_0 > \frac{1}{2}$.
	Notice that there is only one term in function $f_0$,
	then a natural guess for the representation is to consider a branching Brownian motion with exactly one offspring particle at every arrival time. This can be seen as a Brownian motion $W$ equipped with a sequence of random time mark $(T_i)_{i = 1, \cdots, N_T}$,
	where
	$$
		T_i ~:=~ T \wedge \sum_{j=1}^i \tau^{j,1}, ~~~N_T ~:=~ \inf \Big\{i ~: T_i \ge T \Big\}. 
	$$
Notice that for any $t >0$ and bounded measurable function $\phi: \R \to \R$, one has
	$$
		\partial^2_{xx} \E \big[ \phi(x + W_t) \big]
		~=~
		\E \Big[ \phi(x +W_t) \frac{W_t^2 - t}{t^2} \Big].
	$$
	Then arguing as in Theorem \ref{thm:main1}, we may expect that
	$u(0,x_0) = \E[ \widehat \psi \big]$, with
	$$
		\widehat \psi := g(x +W_T) \frac{1}{\Fb(T- T_{N_T-1})} \prod_{i=1}^{N_T - 1}
		c_0 \frac{(W_{T_{i+1}} - W_{T_i})^2 - (T_{i+1} - T_i)}{ (T_{i+1} - T_i)^2 \rho(T_i - T_{i-1})},
	$$
	provided that $\widehat\psi$ is integrable.
	However,	the integrability of $\widehat\psi$ could fail in general. 
	For simplicity, let $g \equiv 1$, and notice that $\Fb \le 1$.
	Then by taking conditional expectation, one has, 
	for some constant $C>0$ and $c_1 := \E \big[ \big| c_0 (W_1^2 - 1) \big| \big]$, that
	\b*
		\E \big[ | \widehat \psi | \big]
		&\ge&
		\E \Big[ \prod_{i=1}^{N_T-1} \frac{c_1}{(T_{i+1} - T_i) \rho(T_i - T_{i-1})} \Big] \\
		&\ge&
		\E \Big[ \frac{c_1}{\rho(T_1)} \frac{c_1}{(T_2- T_1) \rho(T_2 - T_1)} \frac{c_1}{T_3 - T_2} \1_{\{T_1 \le T/2, T_2 - T_1 < T/2, T_3 - T_2 \ge T\}}
		\Big] \\
		&\ge&
		C \E \Big[ \frac{c_1}{(T_2- T_1) \rho(T_2 - T_1)}  \1_{\{T_2 - T_1 < T/2\}} \Big]
		~=~
		C \int_0^{T/2} \frac{1}{t} dt
		~=~
		\infty.
	\e*

	Of course, for linear PDEs as in \eqref{eq:PDE_1}, one can simulate a Brownian motion with volatility coefficient $1 + 2 c_0$ whenever $1 + 2 c_0 > 0$ to obtain the solution.
	But it is not the case for general fully nonlinear PDEs.

\paragraph{On the representation results by BSDE}
	Another probabilistic representation of semilinear parabolic PDE is the Backward Stochastic Differential Equation (BSDE) proposed by Pardoux and Peng \cite{PardouxPeng}.
	Namely, given a classical solution $u$ of semilinear PDE \eqref{eq:PDE}, we define	
	$$
		(Y_t, Z_t)
		~:=~
		\big( u(t,\Xb^{0,x_0}_t), \sigma Du(t, \Xb^{0,x_0}) \big).
	$$
	Then $(Y, Z)$ provides a solution to BSDE
	$$
		Y_t ~=~ g(\Xb^{0,x_0}_T) 
		+ \int_t^T f \big(s, \Xb^{0,x_0}_s, Y_s, \sigma^{-1}(s, \Xb^{0,x_0}_s) Z_s \big) ds 
		- Z_s dW_s,
		~~t \in [0,T],~\P\mbox{-a.s.}
	$$	
	
	Based on the discretization technique on the BSDE, 
	one can then obtain a probabilistic numerical solution for semilinear parabolic PDEs, see e.g. Bouchard and Touzi \cite{BouchardTouzi}, and Zhang \cite{Zhang}, etc.
	Generally speaking, these numerical schemes for BSDE need a (time-consuming) simulation-regression technique
	to compute the conditional expectation appearing in the schemes.

	Our representation result induces a pure Monte Carlo simulation algorithm, 
	which avoids the regression procedure in the numerical schemes of BSDEs.
	Nevertheless, our numerical method provides only the solution of PDE at time $0$,
	and it needs some restrictive conditions on the coefficient functions $f$ such as Assumption \ref{assum:integrability_explicit} to obtain a finite variance estimator.
	We will provide more numerical examples as well as some variance reduction techniques in Section \ref{sec:MC} below.

\section{A Monte Carlo algorithm}
\label{sec:MC}

\subsection{The implementation of the numerical algorithm}

	The above representation result in Theorem \ref{thm:main1}
	induces a Monte Carlo algorithm to compute the solution of PDE \eqref{eq:PDE},
	by simulating the random variable $\psi$ or $\psih$.
	We provide here some discussion on the implementation of the numerical algorithm.

\paragraph{The choice of density function $\rho$}

	As discussed in Remark \ref{rem:integrability_explicit}, to ensure Assumption \ref{assum:integrability_explicit},
	a necessary condition is to choose $\rho(t) \ge C t^{-1/2}$.
	A natural candidate as distribution, which is also easy to be simulated, 
	is the gamma distribution $\Gamma(\kappa, \theta)$, with $\kappa \le \frac{1}{2}$,
		whose density function is given by
		\be \label{eq:gamma_density}
			\rho_0(t) ~=~ \frac{1}{\Gamma(\kappa) \theta^{\kappa}} t^{\kappa -1} \exp(- t/\theta) \1_{\{t > 0\}},
		\ee
		where $\Gamma(\kappa) := \int_0^{\infty} s^{\kappa-1} e^{-s} ds$.
		In particular, one has 
		$$
			\Fb_k 
			~:=~ \int_{\Delta T_k}^{\infty} \rho_0(t) dt 
			~=~ 1- \frac{\gamma(\kappa, \Delta T_k/\theta) }{\Gamma(\kappa)},
			~~~\mbox{where}~~
			\gamma(\kappa, t) := \int_{0}^{t} s^{\kappa -1} e^{-s} ds.
		$$

\paragraph{Complexity} 

	The dimension $d$ of the problem, the choice of $(p_{\ell})_{\ell \in L}$ and $\rho$ will of course influence the complexity of
	algorithm.
	First, the complexity is proportional to the number of particles in the branching process,
	i.e. $\# \Kcb_T$,
	and for each particle, the complexity of simulation and calculation is of order $C d^2$.
	Let us denote $n_0 := \sum_{\ell \in L} p_{\ell} |\ell |$ and  $m(t) := \E \big[ \# \Kcb_t \big]$.

	\begin{Proposition}
		\rmi 
		The function $m(t)$ is given by
		\be \label{eq:m_t_value}
			m(t) ~=~ \sum_{k=0}^{\infty} n_0^k F^{*,k}(t),
			~~~\mbox{where}~~F^{*,k}(t) := \P \big[\tau^{1,1} + \cdots + \tau^{1,k} < t \big].
		\ee
		
		\noindent \rmii Let $\rho$ be given by \eqref{eq:gamma_density}, then
		$F^{*,k}(t) = \frac{1}{\Gamma(k \kappa)} \gamma(k \kappa, x/\theta)$ and hence
		$$
			m(t) ~=~  \sum_{k=0}^{\infty} \frac{\gamma(k \kappa, t/\theta) n_0^k}{\Gamma(k \kappa)}.
		$$
	\end{Proposition}
	\proof \rmi Using Lemma 4.4.3 of Athreya and Ney \cite{AthreyaNey}, one has that
	satisfies the equation $m(t) ~=~ 1 ~+~ n_0 \int_0^t m(t-s) \rho(s) ds$, whose solution is given explicitly by
	\eqref{eq:m_t_value}.
	Further, when $\rho$ is the density function of Gamma distribution, the function $F^{*,k}(t)$ can be computed explicitly.
	\qed

\subsection{A high dimensional numerical example}
	We first focus  on a simple numerical example in high dimension.
	Let $(\mu, \sigma) \equiv (0, \sigma_0)$ for some constant matrix $\sigma_0 = \frac{1}{\sqrt{d}} \I_d$,
	and $ f(t,x,y,z)= k(t,x) + c y (b \cdot z)$,
	where $b := \frac{1}{d} (1+\frac{1}{d}, 1+\frac{2}{d}, \cdots, 2)$ and 
	$$
		k(t,x) 
		:=
		\cos( x_1 +\cdots + x_d) 
		\Big( \alpha + \frac{\sigma^2}{2}  + c \sin(x_1 +\cdots + x_d)  \frac{3d+1}{2d} e^{ \alpha (T-t)} \Big) 
		e^{ \alpha (T-t)}.	
	$$
	With terminal condition $g(x) = \cos( x_1+ \cdots + x_d)$,
	the explicit solution of semilinear PDE \eqref{eq:PDE} is given by
	$$
		u(t,x) =  \cos(x_1+ \cdots + x_d) e^{\alpha (T-t)}.
	$$
	In our numerical experiment, we set $\alpha=0.2$, $c= 0.15$, $T=1$, and $x_0 =0.5 \un_d$,
	where  $\un_d$ stands for  the unit vector in $\R^d$ for $d=5,10$ and $20$. 
	We would like to emphasize that, to the best of our knowledge, no alternative methods are available for  solving  such a high-dimensional semilinear PDE. 
	In Table \ref{tab:Comp}, we report the analytic solution of the semilinear PDE
	and that of the corresponding linear PDE by setting $c = 0$.
	The different results indicate that the nonlinearity  term has an impact. 
\begin{table}[h!]
  \begin{center}
    \begin{tabular}{|l|l|l|l|}
      \hline
      Dimension & 5 &  10 &  20      \\ \hline
      Linear Solution  &  -1.0436  & 0.3106  &   -0.9661   \\ \hline
      Non linear solution &  -0.97851 &  0.34646 & -1.0248  \\ \hline
    \end{tabular}
  \end{center}
  \caption{ Analytical solution for the linear PDE (i.e., $c=0$) versus analytical solution for the semilinear PDE in $d=5,10$ and $20$.}
   \label{tab:Comp}
\end{table}

	For numerical implementations, we use gamma distribution \eqref{eq:gamma_density}, with $\kappa = 0.5$ and $\theta=2.5$.
	On each test performed, a computation is achieved with $n$ particles. An estimation $E$ with $n$ particles is then calculated. The standard deviation  of $E$ is estimated with $1000$ runs of $n$ particles and its log-plot is reported below on the different figures for different values of $n$. We also  show on some figures the convergence of the solution obtained on the average of the 1000 runs. 

\begin{figure}[h!]
  \centering
  \includegraphics[width=6cm]{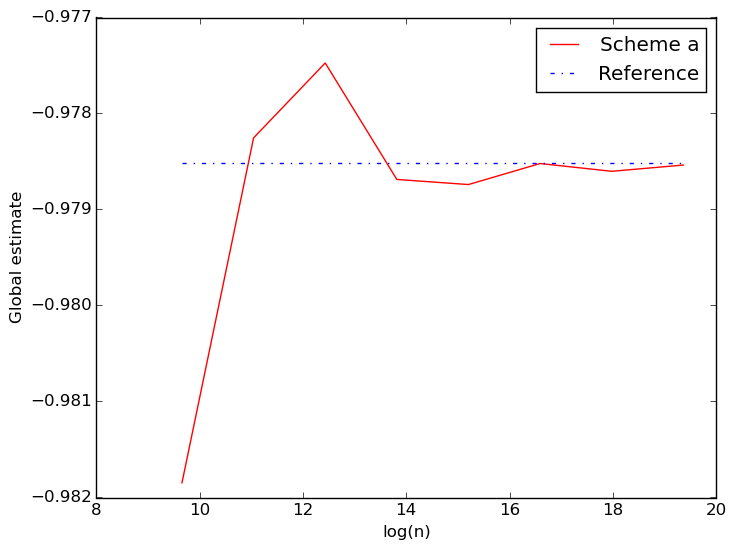}
  \includegraphics[width=6cm]{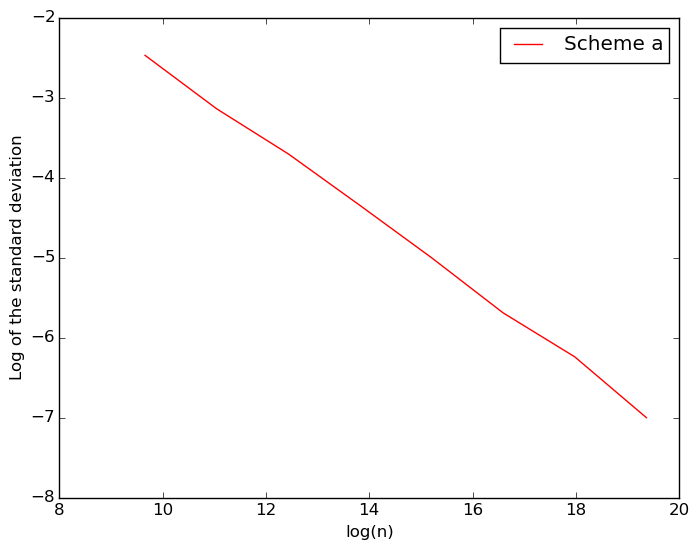}
  \caption{Estimation and standard deviation observed in $d=5$ depending on the log of the number of particles used.}
  \label{FigCaseHD1}
\end{figure}

\begin{figure}[h!]
  \centering
  \includegraphics[width=6cm]{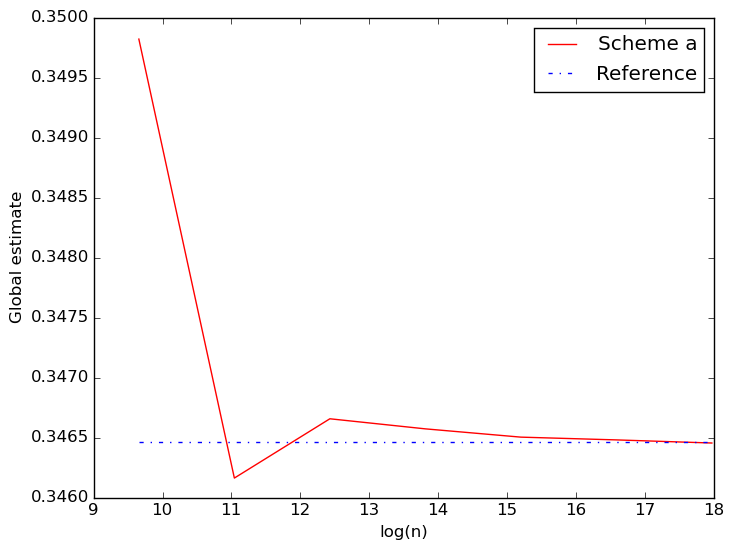}
  \includegraphics[width=6cm]{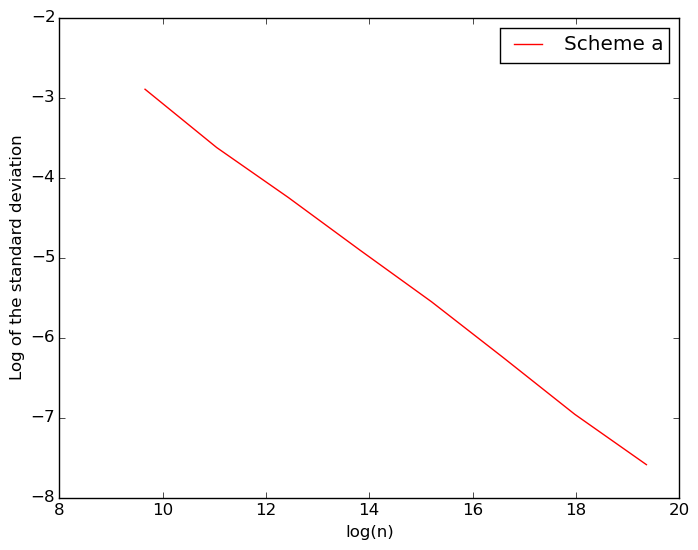}
  \caption{Estimation and standard deviation observed in $d=10$ depending on the log of the number of particles used.}
  \label{FigCaseHD2}
\end{figure}

\begin{figure}[h!]
  \centering
  \includegraphics[width=6cm]{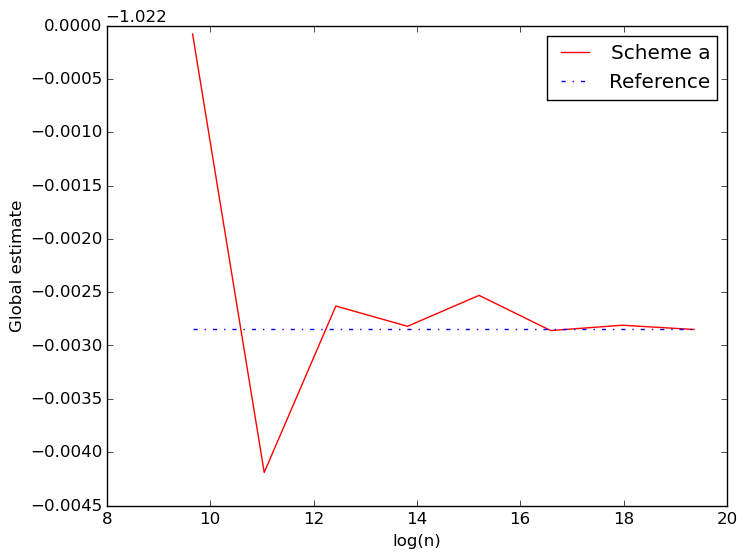}
  \includegraphics[width=6cm]{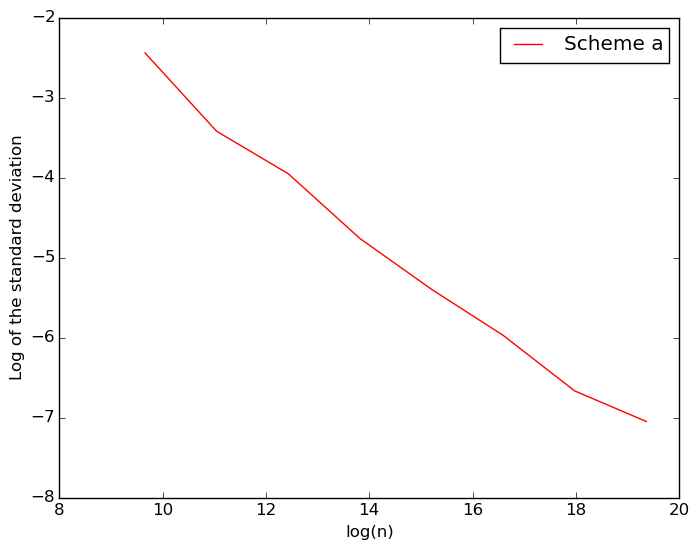}
  \caption{Estimation and standard deviation in $d=20$ depending on the log of the number of particles used.}
  \label{FigCaseHD3}
\end{figure}

 On Figures \ref{FigCaseHD1}, \ref{FigCaseHD2}, \ref{FigCaseHD3}, we illustrate that the Monte  Carlo method converges easily to our analytic solution. Computational costs are estimated on one core of  a Laptop core I7 processor 2.2 GHz and are reported in Table \ref{tab:table1} for a number of simulations equal to $96000$ permitting 
to get a solution with an error less than $0.1 \%$.

\begin{table}[h!]
  \begin{center}
    \begin{tabular}{|l|l|l|l|}
      \hline
      Dimension &  5  & 10   &  20 \\ \hline
      Time      &  550   &  717     & 956  \\ \hline
    \end{tabular}
  \end{center}
   \caption{Computational time in seconds for $96000$ trajectories computed 1000 times on one core for $\kappa =0.5$, $\theta=2.5$.}
    \label{tab:table1}
\end{table}

\section{Some extensive tests}
\label{sec:Simu2}

	This section is devoted to additional  tests. Having illustrated previously that our algorithm is efficient  for solving high-dimensional semilinear PDEs, we focus on some examples from dimension $1$ to $3$. 
	Note that our results have been benchmarked against a Finite Difference method in $d=1$ and $d=2$.  Unfortunately,  the finite difference method is no more available in $d=3$. All our numerical examples  share the following characteristics: 
	$\mu(t,x)= 1- x$, $\sigma \equiv 0.5 \I_d$ and $x_0 = \un$.
	$T$ is chosen equal to $1$,
	$g(x) = (\frac{1}{d} \sum_{i}^d x(i)- 1)^{+}$ for $x \in \R^d$.
	Notice that with the above coefficients, SDE \eqref{eq:SDE_Xb} is a linear SDE, whose solution can exactly simulated:
	\be \label{eq:SDE_explicite}
		\Xb_t^{0,x} = (1-e^{- t}) \un + e^{-t} x + \sigma \sqrt{\frac{1 -e^{-2t}}{2}} Z,
		~~~&&Z \sim N(0,\I_d).
	\ee
	The Malliavin weight  used in the algorithm can be computed explicitly and is given by $ Z / \big(\sigma \sqrt{\frac{e^{2 t}-1}{2}}\big)$.
	We will compare numerical results from four different schemes.
	\begin{itemize}
		\item ({\bf  scheme  a}) using the representation \eqref{eq:def_psi} with the explicit solution \eqref{eq:SDE_explicite} of the SDE \eqref{eq:SDE_Xb}. 
		\item ({\bf scheme b}) using the representation \eqref{eq:def_psi_exact} with freezing coefficient techniques.
		\item ({\bf scheme c})   using the representation \eqref{eq:def_psi_exact},  enhanced by the resampling scheme (see Appendix for more details). 
		\item ({\bf scheme d}) using the representation \eqref{eq:def_psi}, enhanced by the resampling scheme.
	\end{itemize}

	The density function $\rho$ is that of the gamma law with parameters $\kappa$ and $\theta$. 
	If not indicated, the parameters of the law are set to $\kappa =0.5$ and $\theta=2.5$ and  the probability $p_l$ are chosen equal.
	On each test, a calculation is achieved with $n$ particles (starting with $n=1562$ for scheme {\bf a} and  with $n=100000$ for schemes {\bf b} and {\bf c}). The $n$ particles are
shared on 96 processors and each processor $i$ calculates an estimation $E_i$  of the solution  with $\frac{n}{96}$ particles. Then an estimation $E$ with $n$ particles is achieved with
$E = \frac{1}{96} \sum_i E_i$. When importance sampling is used, in order to avoid communications that breaks parallelism, it is used on each processor so with $\frac{n}{96}$ particles on each processor.
The standard deviation  of $E$ is estimated with 1000 runs of $n$ particles and its log is reported on the different figures below for different values of $n$. We expect that  by quadrupling the values of $n$, the standard deviation  $std$ divides by a factor $2$ and the plot $(\log(n),\log(std))$ should be linear with a slope equal to $-\frac{1}{2}$. The theoretical rate of convergence is also plotted on each figure (as in our previous example, the solutions are obtained on the average of the 1000 runs).

\subsection{Some examples in one space dimension}

\begin{figure}[h!]
  \centering
  \includegraphics[width=6cm]{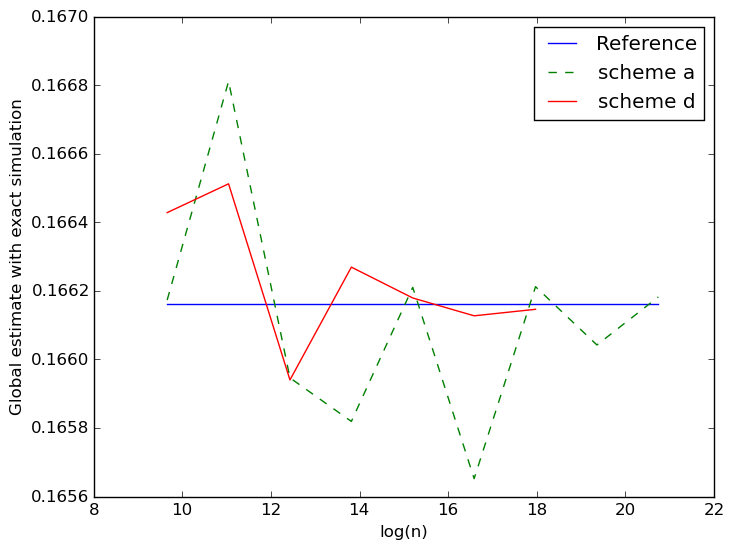}
  \includegraphics[width=6cm]{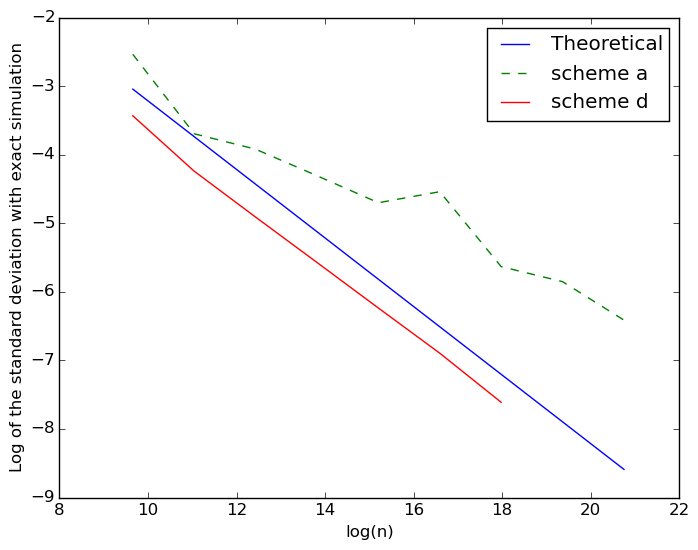}
  \caption{Estimation and standard deviation obtained  in $d=1$ for  $f(t,x,y,z) =  0.2 z^2$ with and with out importance sampling.}
  \label{FigCase1D41}
\end{figure}

\begin{itemize}
\item For $d=1$, we take   $ f(t,x,y,z) :=  0.2 y^2 + 0.3 y^3$.
  Results on Figure \ref{FigCase1D1} show that the method converges. Scheme {\bf a} is far more effective than scheme {\bf b} and that the importance sampling of scheme {\bf c} is effective. The log of the standard deviation decreases for all schemes  linearly with  the log of the particle number as  predicted by the theory.
Note that the computational cost for 1000 runs with 25000 particles on one core is equal to  490 seconds for scheme {\bf a}, 200 seconds with scheme {\bf b} and  260 seconds with scheme {\bf c}.
  
  \begin{figure}[h!]
    \centering
    \includegraphics[width=6cm]{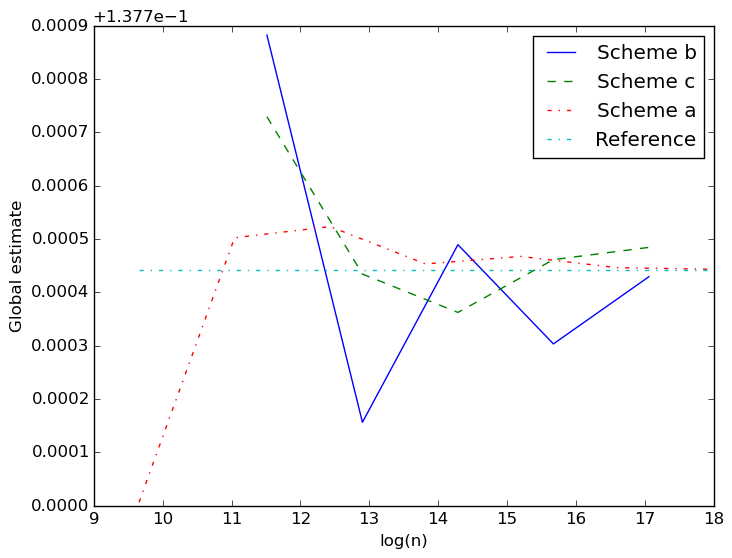}
    \includegraphics[width=6cm]{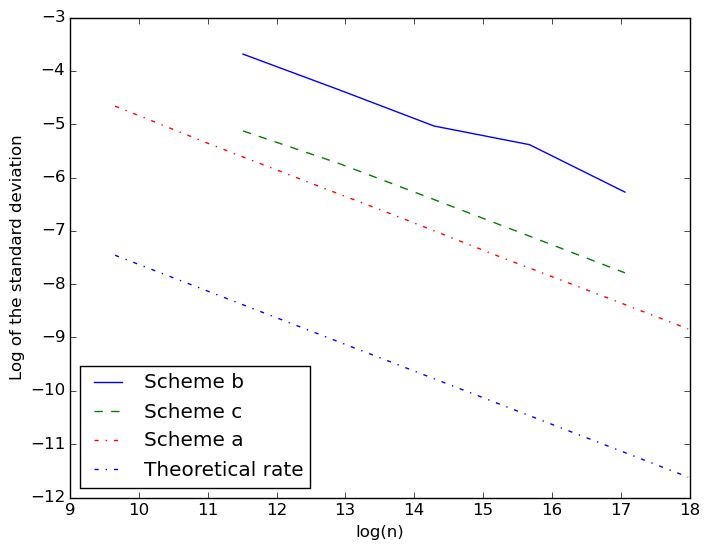}
    \caption{Estimation and standard deviation obtained in $d=1$ for $ f(t,x,y,z) :=  0.2 y^2 + 0.3 y^3$}
    \label{FigCase1D1}
  \end{figure}
\item As a second example in $d=1$, we take  a Burgers type nonlinearity  $ f(t,x,y,z) =  0.15 yz$.
  Results on Figure  \ref{FigCase1D2} show that all the schemes converge to our numerical finite difference solution. Note that the computational cost for 1000 runs with 25000 particles on one core is roughly equal to 200 seconds for scheme {\bf a}, 100 seconds for scheme {\bf b}, 300 seconds for scheme {\bf c}. 
 \begin{figure}[h!]
    \centering
    \includegraphics[width=6cm]{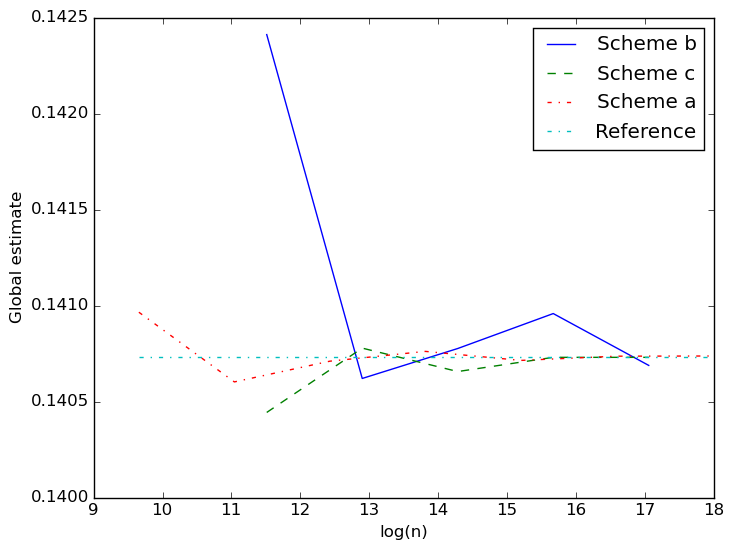}
    \includegraphics[width=6cm]{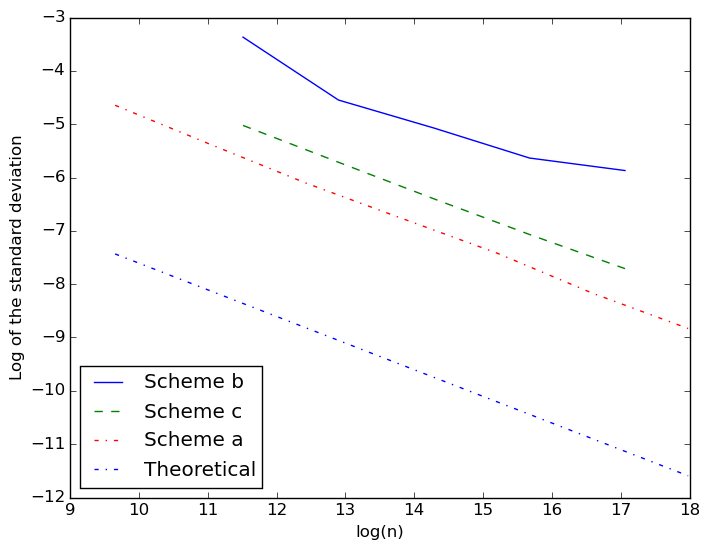}
    \caption{Estimation and standard deviation obtained in $d=1$ for $ f(t,x,y,z) =  0.15 yz$.}
    \label{FigCase1D2}
  \end{figure}
\item  As a third example in $d=1$, we keep the same nonlinearity with $ f(t,x,y,z) =  0.3 yz$. We expect that the variance of the results will be higher than in the previous case. This  is observed in Figure \ref{FigCase1D5}.
  Scheme {\bf a}  still converges. Scheme {\bf b}  converges slowly and Importance Sampling of scheme {\bf c} permits to get faster convergence and to recover the good rate in the log of the standard deviation decay. The computational times are the same as in our previous test.
  \begin{figure}[h!]
    \centering
    \includegraphics[width=6cm]{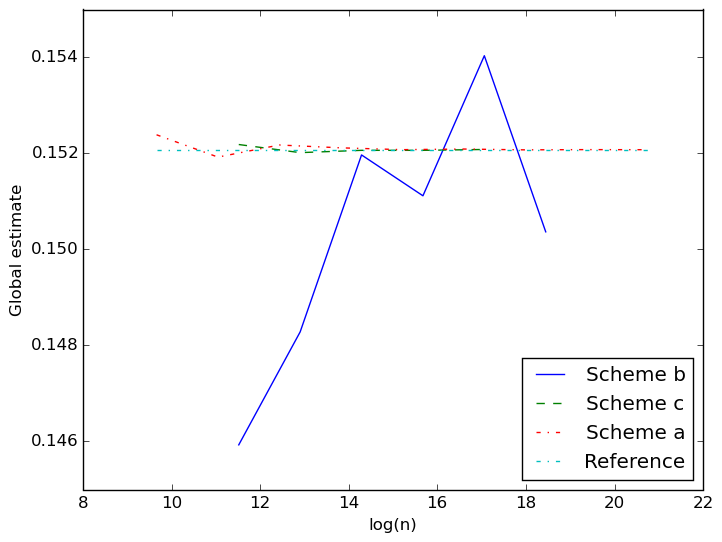}
    \includegraphics[width=6cm]{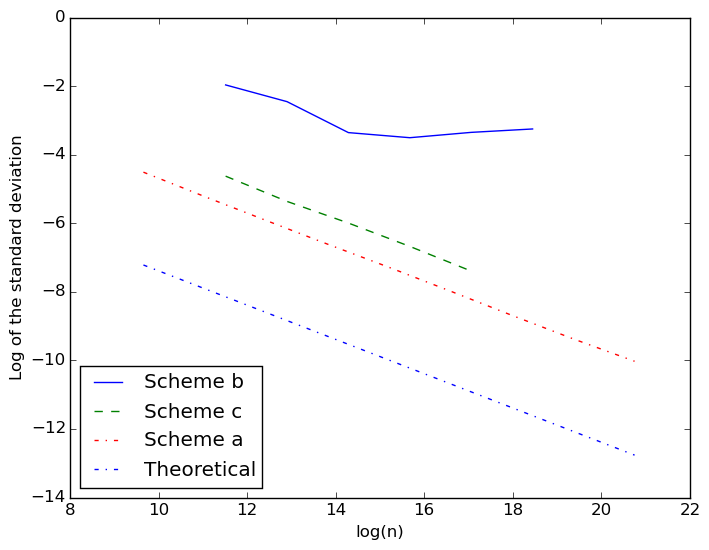}
    \caption{Estimation and standard deviation obtained in $d=1$ for $ f(t,x,y,z) =  0.3 yz$. }
    \label{FigCase1D5}
  \end{figure}
\item As a fourth example in $d=1$, we take a nonlinearity with $ f(t,x,y,z) =  0.08 z^2$. Results are shown in Figure \ref{FigCase1D7}. The importance sampling of scheme {\bf c} is required to achieve proper convergence. Scheme {\bf a} converges quickly.  Computational times are the same as before (same type of branching).
  \begin{figure}[h!]
    \centering
    \includegraphics[width=6cm]{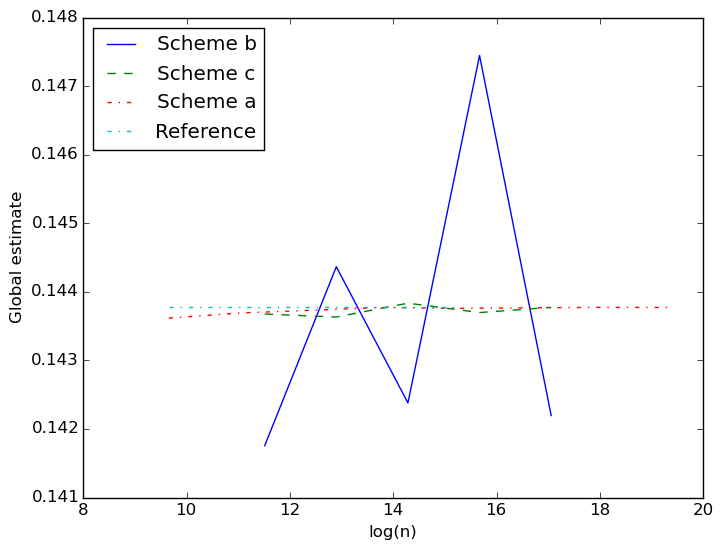}
    \includegraphics[width=6cm]{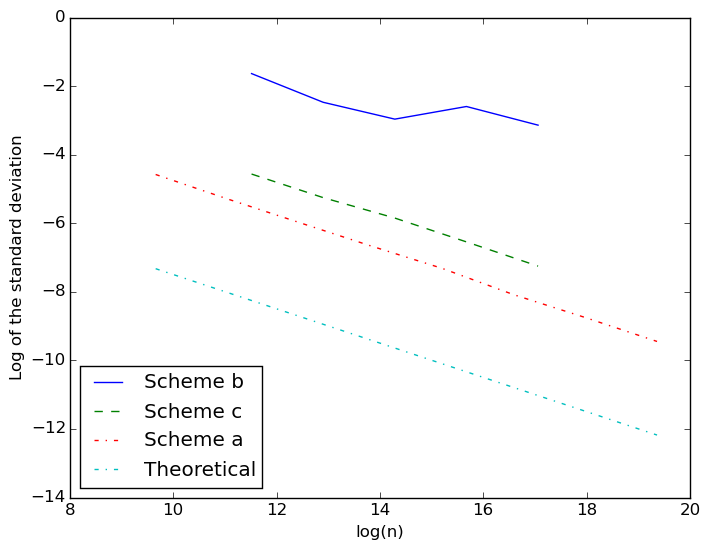}
    \caption{Estimation and standard deviation obtained in $d=1$ for $ f(t,x,y,z) =  0.08 z^2$. }
    \label{FigCase1D7}
  \end{figure}
\item As a last example in $d=1$, we keep the same type of nonlinearity  $ f(t,x,y,z) =  0.2 z^2$. Schemes {\bf b} and {\bf c} don't converge anymore.
  We only test scheme {\bf a} using different values for the parameters $\kappa$ and $\theta$ (see Figure \ref{FigCase1D4}). The change in  $\theta$ does not seem to change convergence properties.  The change in $\kappa$ (from $0.5$ to $0.4$) does not seem to modify our results. However, some tests, not reported here, show that the variance can increase a lot using  $\kappa$ around $0.25$.
  \begin{figure}[h!]
    \centering
    \includegraphics[width=6cm]{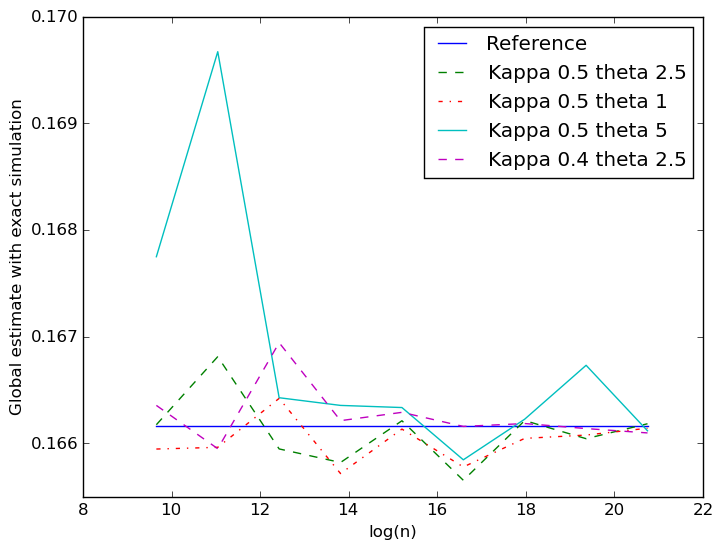}
    \includegraphics[width=6cm]{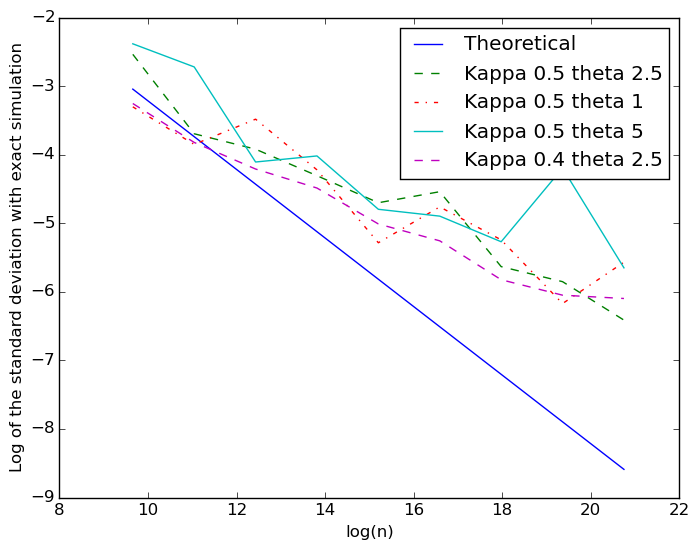}
    \caption{Estimation and standard deviation obtained in $d=1$ for $ f(t,x,y,z) =  0.2 z^2$. Different values for $\kappa$ and $\theta$ are used. }
    \label{FigCase1D4}
  \end{figure}
  Then,  as the average jump size is proportional to $\theta$, it is more efficient to take some quite high values for $\theta$ in order to reduce the computational time. For the same reason, it is optimal to choose a $\kappa$ equal to $0.5$.
  In Table \ref{table2}, we report the computational time, associated to different choices of  $(\kappa, \theta)$, as a  multiplicative factor of the computational effort  with benchmark parameters $\kappa=0.5$, $\theta=2.5$.
  \begin{table}[h!]
    \caption{Computational time, associated to different choices of  $(\kappa, \theta)$, as a  multiplicative factor of the computational effort  with benchmark parameters $\kappa=0.5$, $\theta=2.5$}
    \label{table2}
    \begin{center}
      \begin{tabular}{|l|l|l|l|l|l|}
        \hline
        $\kappa$ &  0.5 & 0.5 & 0.5 & 0.4 & 0.4 \\ \hline
        $\theta$ &  1 & 2.5 & 5 & 2.5 & 5.\\ \hline
        Time     &  6.63  & 1  & 0.49 & 2.85 & 1.02 \\ \hline
      \end{tabular}
    \end{center}
  \end{table} 
\end{itemize}
We notice that for all the parameters, the decay in the variance is far from the expected theoretical one (see Figure \ref{FigCase1D4}). We then use our benchmark parameters
and compare the results obtained using scheme {\bf a} and scheme {\bf d} (importance sampling is used here). Results are reported on Figure \ref{FigCase1D41}. They illustrate that the importance sampling method allows to improve the convergence rate.

\subsection{Some examples in two space dimensions}

Although the efficiency of our algorithm was illustrated on our previous experiments, this Monte-Carlo method cannot compete a PDE  deterministic methods in $d=1$.  In this section, we focus on $d=2$, where  advantages of PDE implementation remain but are not so obvious.
\begin{itemize}
\item For the first example in $d=2$, we take $ f(t,x,y,z) := 0.15 y \un .z$.  On Figure  \ref{FigCase2D2}, we give the results obtained using our three schemes showing that Importance Sampling is needed.  Note that  the computation cost for 1000 runs with 25000 particles on one core is roughly equal to 230 seconds for scheme {\bf a}, 90 seconds for scheme {\bf b}, 580 seconds for scheme {\bf c}.

  \begin{figure}[h!]
    \centering
    \includegraphics[width=6cm]{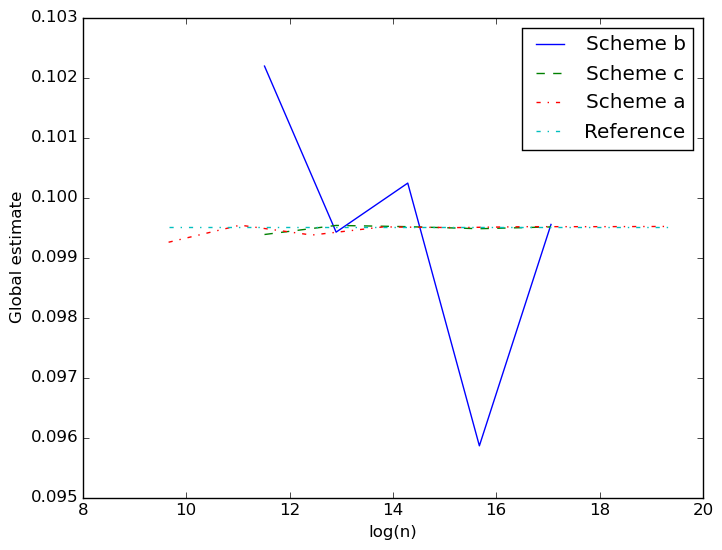}
    \includegraphics[width=6cm]{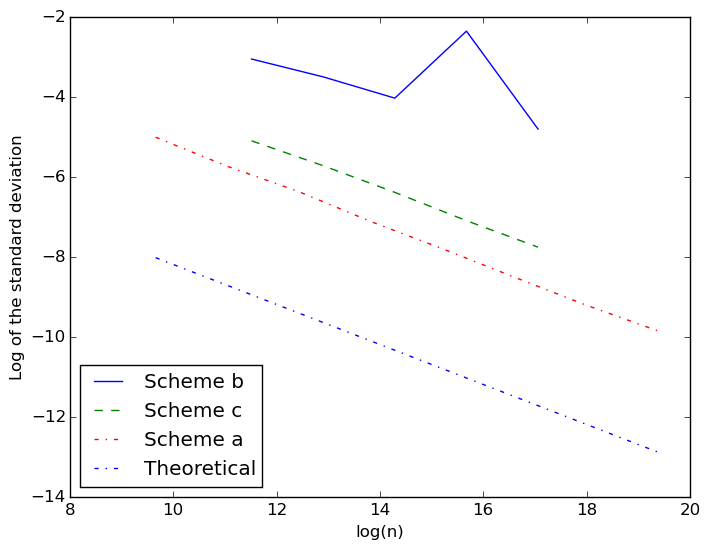}
    \caption{Estimation and standard deviation obtained in $d=2$ for $ f(t,x,y,z) := 0.15 y \un .z$.}
    \label{FigCase2D2}
  \end{figure}
\item For the second example in $d=2$, we take $ f(t,x,y,z) := 0.04 (z.\un)^2$.
  The convergence of Scheme {\bf a} is easily achieved while Scheme {\bf b} converge poorly as shown in Figure \ref{FigCase2D3}. Importance sampling method improve the convergence. Computational costs are the same as in our first $d=2$ tests.
  \begin{figure}[h!]
    \centering
    \includegraphics[width=6cm]{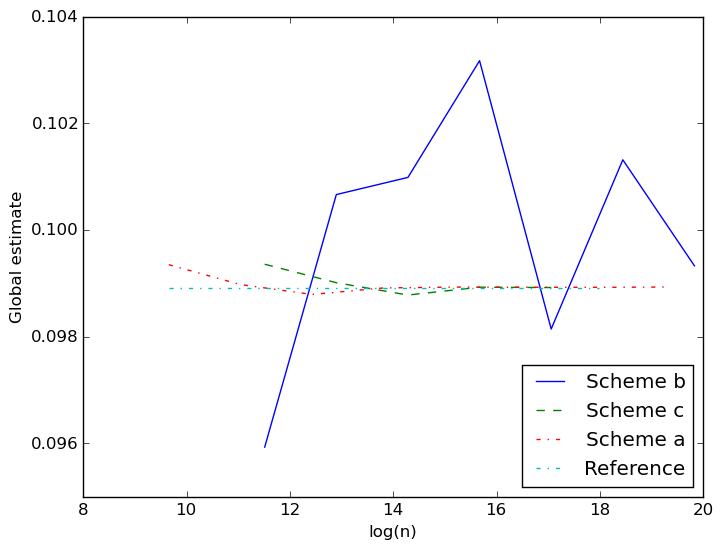}
    \includegraphics[width=6cm]{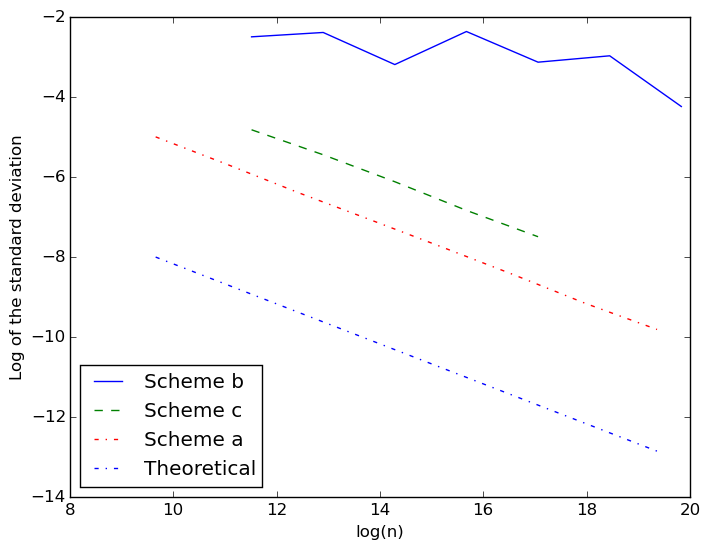}
    \caption{Estimation and standard deviation obtained in $d=2$ for $ f(t,x,y,z) := 0.04 (z.\un)^2$.}
    \label{FigCase2D3}
  \end{figure}
\item For the third example, we test the influence of the coefficients  on Scheme {\bf a} for a non linearity $ f(t,x,y,z) := K (z.\un)^2$ with $K=0.05$, $K=0.1$, $K=0.2$. Using Scheme {\bf b} and { \bf c}, we cannot get proper convergence due to high variances observed. On Figure \ref{FigCase2D4_1}, we give the convergence obtained with the different $K$ values and on Figure \ref{FigCase2D4_2} the standard deviation associated.
  As the coefficients grow, the variance of the results gets higher preventing the method from converge  when $K=0.2$.
  \begin{figure}[h!]
    \centering
      \includegraphics[height=1.2in]{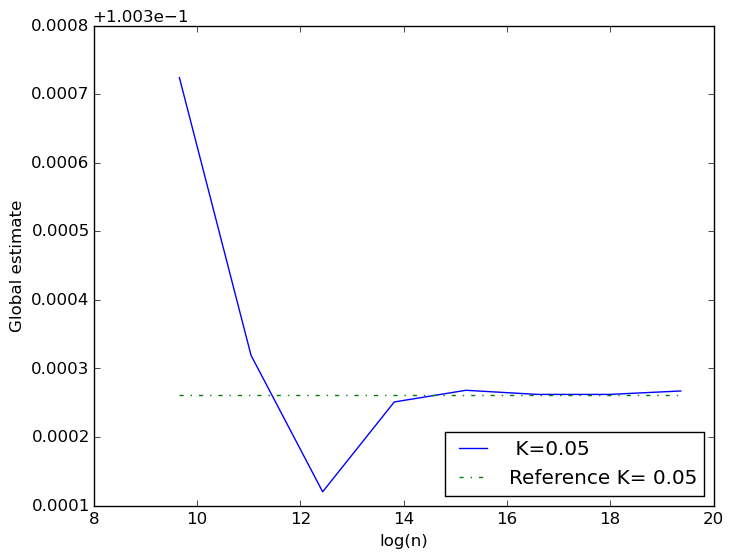}
      \includegraphics[height=1.2in]{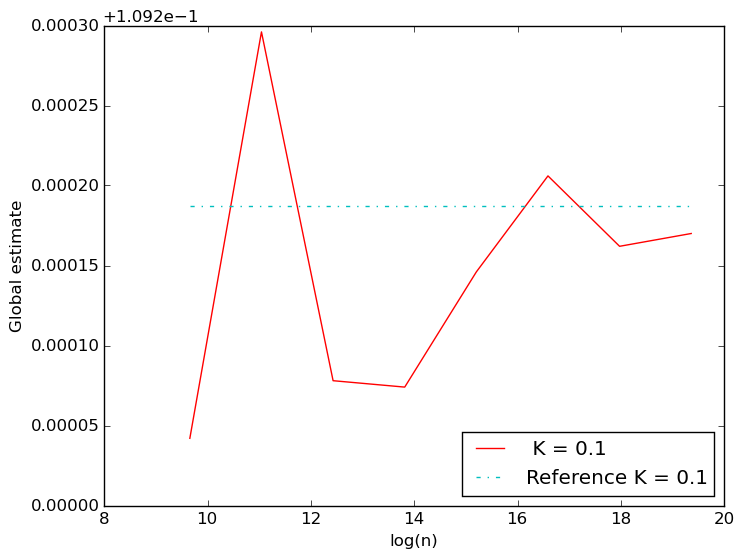}
      \includegraphics[height=1.2in]{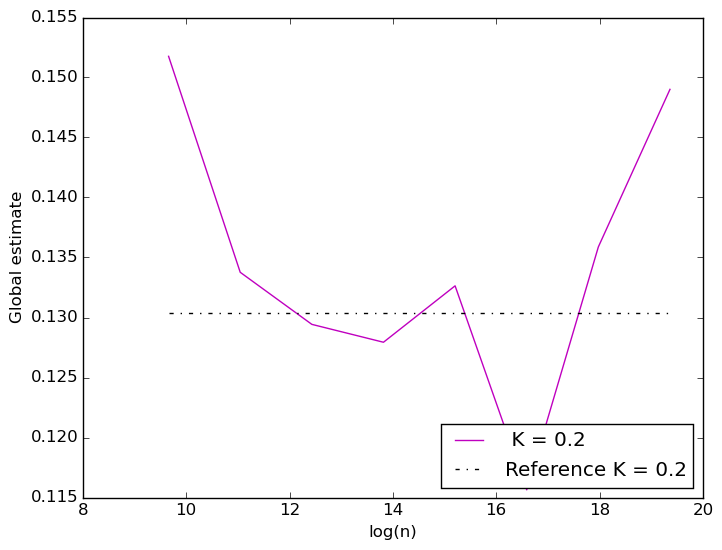}
    \caption{Convergence of scheme {\bf a} for different $K$ values.}
    \label{FigCase2D4_1}
  \end{figure}

  \begin{figure}[h!]
    \centering
      \includegraphics[height=1.2in]{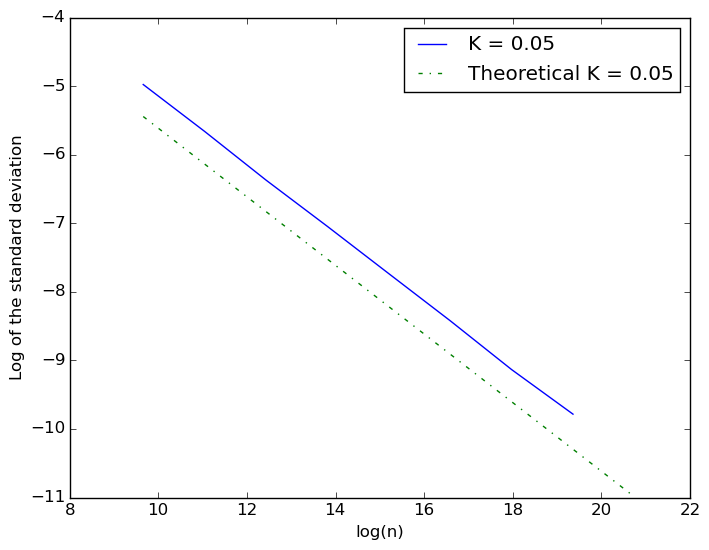}
      \includegraphics[height=1.2in]{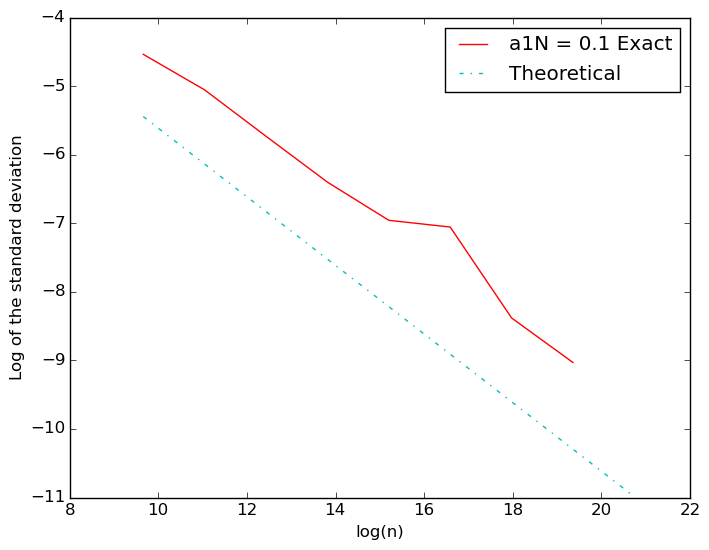}
      \includegraphics[height=1.2in]{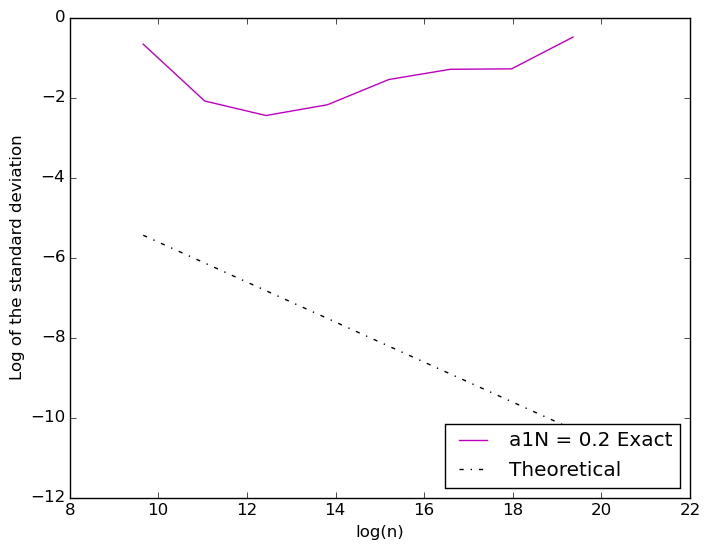}
    \caption{Standard deviation of the scheme {\bf a} for different $K$ values.}
    \label{FigCase2D4_2}
  \end{figure}
\item At last we test the influence of the function $g$. The representation of the solution involves the product of $g$ functions so we expect that the variance
  of the result is highly sensitive to the scaling of this function. Here we choose to keep $ f(t,x,r,p) := 0.05 (Du.\un)^2$  and take  different values for the $g$ function.
  On figure \ref{FigCase2D5} we take $g(x) = 2 (\frac{1}{d} \sum_{i}^d x(i)- 1)^{+}$ and  give the convergence of schemes {\bf a} and {\bf b}  and the standard deviation associated.
  Comparing to figure \ref{FigCase2D4_2} ($K=0.05$), we see a net increase in the variance of the result for scheme {\bf a}. When importance sampling is used (scheme {\bf d}) the decay in term of variance is more regular.
  Increasing the function $g$ such that $g(x) = 3 (\frac{1}{d} \sum_{i}^d x(i)- 1)^{+}$, we give the results obtained on figure \ref{FigCase2D6}. Here importance sampling
is really necessary to recover a good rate of convergence.
  \begin{figure}[h!]
    \centering
    \includegraphics[width=6cm]{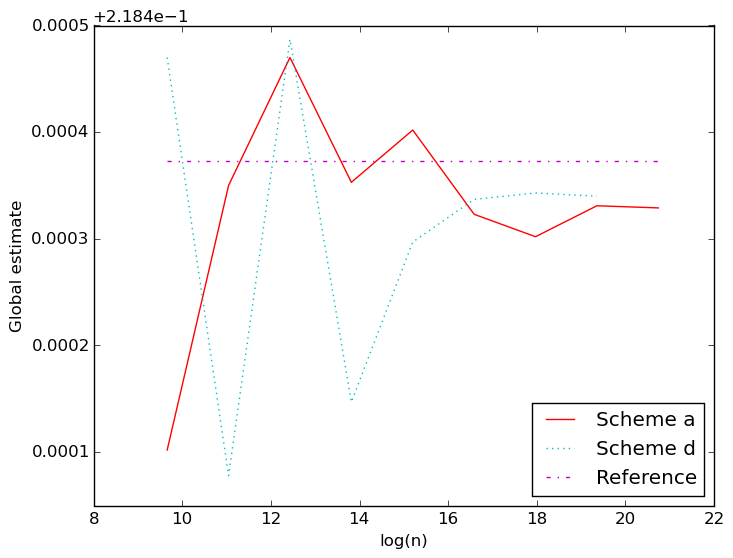}
    \includegraphics[width=6cm]{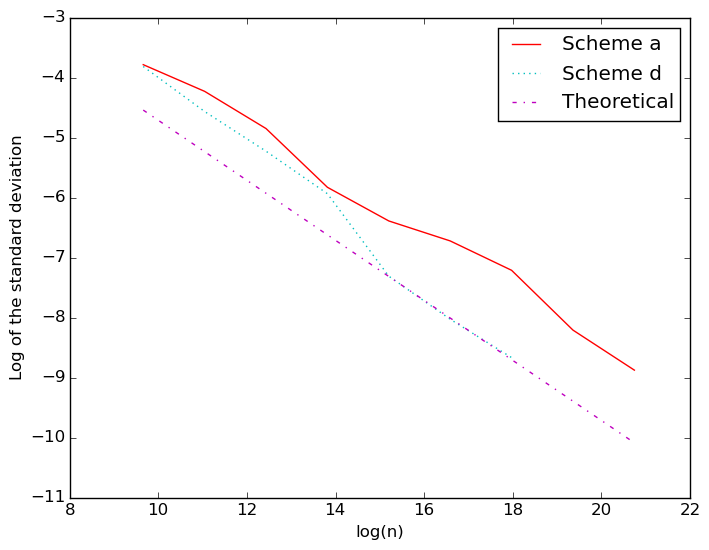}
    \caption{Estimation and standard deviation observed in dimension 2  for  case 4,  $g(x) = 2 (\frac{1}{d} \sum_{i}^d x(i)- 1)^{+}$.}
    \label{FigCase2D5}
  \end{figure}
  \begin{figure}[h!]
    \centering
    \includegraphics[width=6cm]{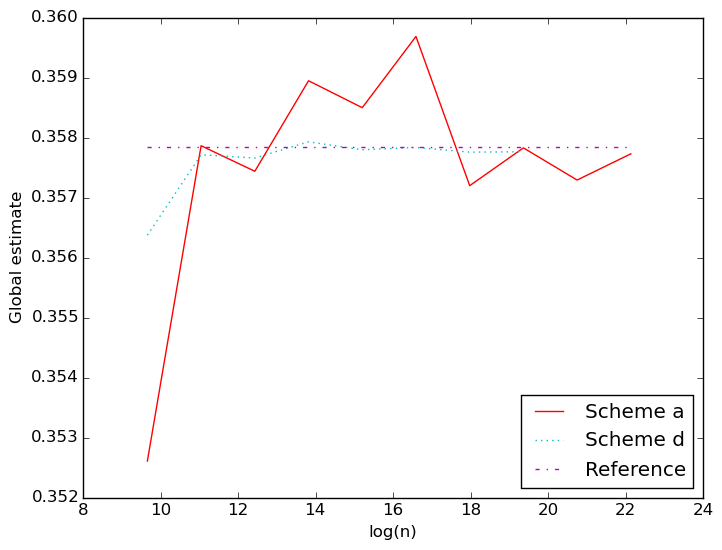}
    \includegraphics[width=6cm]{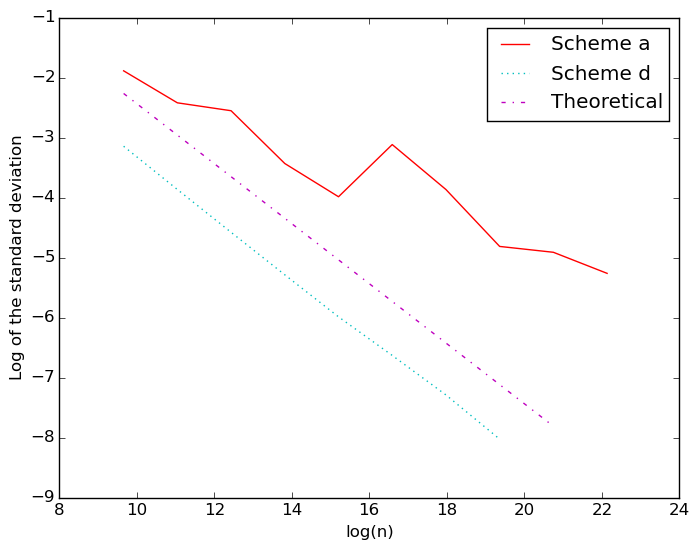}
    \caption{Estimation and standard deviation observed in dimension 2  for  case 4,  $g(x) = 3 (\frac{1}{d} \sum_{i}^d x(i)- 1)^{+}$.}
    \label{FigCase2D6}
  \end{figure}
\end{itemize}

\subsection{An example in three space dimensions}

\begin{figure}[h!]
  \centering
  \includegraphics[width=6cm]{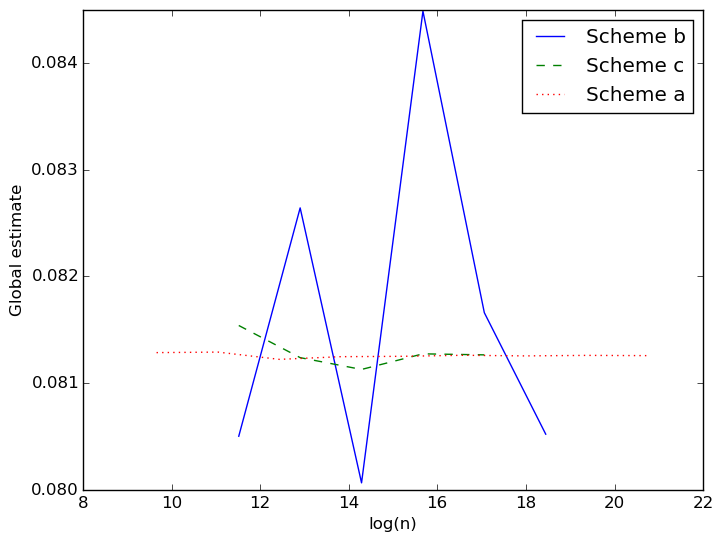}
  \includegraphics[width=6cm]{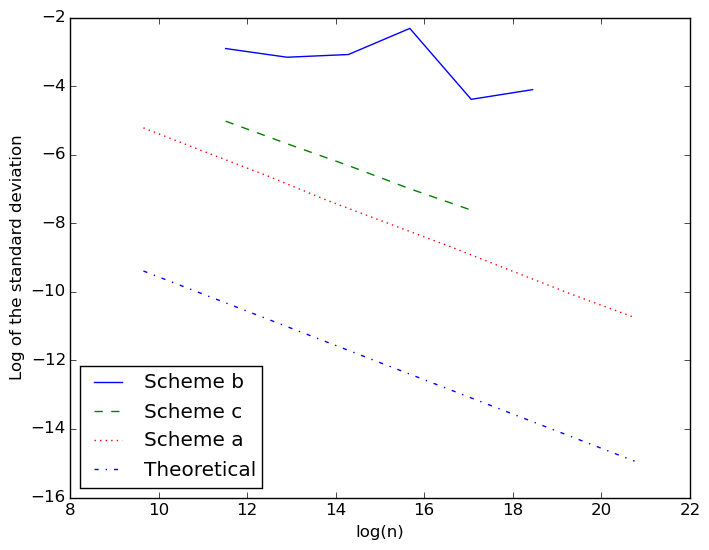}
  \caption{Estimation and standard deviation observed in dimension 3  for  case 1.}
  \label{FigCase3D1}
\end{figure}

	We take $f(t,x,y,z) := 0.15 (z.\un)^2$. Results are given on Figure \ref{FigCase3D1}, still showing that importance sampling is necessary while using discretization of the scheme and that the exact scheme has  a lower variance.

\appendix 

\section{Resampling scheme for branching processes}

	Notice that our estimator \eqref{eq:def_psi} and \eqref{eq:def_psi_exact} are provided as a product of some random variables.
	Then similar to Doumbia, Oudjane and Warin \cite{DOW}, 
	one can use the resampling scheme (or interacting particle systems), see Del Moral \cite{DelMoral}.
	Intuitively, this scheme replaces the expectation of a product by a product of expectations, which potentially stabilizes the Monte-Carlo estimator.

	Let us first introduce the Markov chain $(\Xc _n)_{n\geq 1}$, taking values in 
	$\cup_{p\geq 1} ([0,T]^2\times \R^{2d}\times \{0,\cdots ,m\}\times L)^p$ 
	such that $\Xc_1 =(0,0,x_0,x_0,0,I_0)$ 
	with $I_0=(1,0,\cdots 0)\in \mathbb{N}^{m+1}$ 
	and  for any $n\geq 1$, one defines
	$$
		\Xc _{n+1}
		~:=~
		\Big (\Xc_n, (T_{k^-},T_k,X^k_{T_{k^-}},X^k_{T_k},\theta_k,I_k)_{k\in \cup_{p=1}^n \Kcb^p_T}\Big ).
	$$
	Notice that this Markov chain has an absorbing state since for any $\omega \in\Omega$ there is a generation $n(\omega)$ for which all branches have died (either having no offspring before reaching $T$ or having reached $T$) implying $\Kcb^{n+1}_T(\omega)=\emptyset$ and consequently $\Xc _{n+1}(\omega)=\Xc _n(\omega)$. 
	Then $(\Xc_n)_{n\geq 0}$ is a Markov chain. 
	We next introduce
	\be \label{eq:G}
		G_n( \Xc_n) ~:=
		\Big[\!\! \prod_{k \in  \Kc^n_T} \!\!
		\frac{g(X^k_T) - g(X^k_{T_{k-}}) \1_{\{\theta_k \neq 0\}} }{ \Fb(\Delta T_k)}  \Wc_k \Big ]
		\Big [\!\!\!\! \prod_{k \in  (\Kcb^n_T \setminus \Kc^n_T)}\!\!\!
		\frac{c_{I_k}(T_k, X^k_{T_k})}{p_{I_k}} \frac{\Wc_k}{\rho(\Delta T_k) } \Big ],
	\ee
	so that 
	$$
		\psi ~=~ \prod_{n=1}^{\infty} G_n(\Xc_n).
	$$
	Notice that the above representation consists of a product from contributions from each generation 
	$n \ge 1$.
	Since the number of  generation prior to the maturity $T$ is finite a.s., the last product only involves finite number of terms, a.s.
	We also observe that except for the trivial case of constant function $g$, $\E[|G_n(\Xc_n)|] \neq 0$. 
	By iteration, it is easy to see that
	$$
		\E^{\P_0}[\psi] ~=~ 
		\Big( \prod_{n=1}^{\infty} \E^{\P_{n-1}} \big[\big|G_n(\Xc_n)\big|\big] \Big)
		~\E^{\P_\infty} \Big[ \prod_{n=1}^{\infty} \mbox{sgn}(G_n(\Xc_n)) \Big],
	$$
	where given $\P_0$, one defines $\P_n$ by $\frac{d \P_n}{d \P_{n-1}} := \frac{|G_n(\Xc_n)|}{\E^{\P_{n-1}}[|G_n(\Xc_n)|]}$, for $n \ge 1$.

	The particle algorithm  consists in simulating the dynamics of an interacting particle system of size $N$, 
	$(\xi^{1,N}_p,\cdots \xi^{N,N}_p)$, 
	on $\cup_{p\geq 1} ([0,T]^2\times \R^{2d}\times \{0,\cdots ,m\}\times L)^p$, 
	from step $n=1$ to $n= \infty$ 
	and then to approximate each expectation $ \E^{\P_{n-1}} \big[\big|G_n(\Xc_n)\big|\big]$ by the empirical mean value of the simulation.
	The algorithm can be given as an iteration of the following two steps, initiated by $n=1$,
	\begin{description}
		\item [Selection step]  Given $N$ copies of simulation $(\xi^{i,N}_n)_{i=1, \cdots, N}$ of $\Xc_n$,
		one draws  randomly and independently $N$ particles among the current particle system 
		with a probability $\frac{|G_p(\xi^{i,N}_p)|}{\sum_{j=1}^N |G_p(\xi^{i,N}_p|)} $;

		\item [Evolution step] Each new selected particle evolves randomly and independently according  to the transition of the Markov chain $(\Xc _n)$ between  $n$ and $n+1$.
	\end{description}

	\noindent Finally $u(0,x_0)$ is approximated as a product of empirical averages: 
	\be \label{eq:GammaApprox}
		\prod_{n=1}^{\infty} \Big (\frac{1}{N}\sum_{i=1}^N |G_n(\xi^{i,N}_n) | \Big )
		\Big(\frac{1}{N}\sum_{i=1}^N  \prod_{n=1}^{\infty} \mbox{sgn}\big(G_n(\xi^{i,N}_n) \big) \Big ).
	\ee
	Notice again that, for every simulation $(\xi^{i,N})$, the maturity $T$ is attained for some finite generation,
	then the above product $\prod_{n=1}^{\infty}$ can be restricted to the a finite product $\prod_{n=1}^{n_N}$,
	where $n_N:=\inf \{n \geq 1 \,\vert\,\xi^{i,N}_n  \ \textrm{has reached $T$ for all}\ i=1,\cdots N\}.$

\end{document}